\documentclass[12pt,psamsfonts]{amsart}

\usepackage{amsmath,amssymb}
\usepackage{graphicx}
\usepackage{color}
\usepackage{psfrag} % See the section on figures.
\usepackage{overpic} % See the section on figures.

\newtheorem{theorem}{Theorem}

\newtheorem{example}[theorem]{Example}

\def\R{\mathbb{R}}

\title[]{
Limit Cycles from Planar Piecewise Linear Hamiltonian differential Systems with Two or Three Zones}
\author[C. Pessoa and R. Ribeiro]{}


\begin{document}
 \maketitle

\centerline{\scshape  Claudio Pessoa and Ronisio Ribeiro}
\medskip

{\footnotesize \centerline{Universidade Estadual Paulista (UNESP), Instituto de Bioci\^encias Letras e Ci\^encias Exatas,} \centerline{R. Cristov\~ao Colombo, 2265, 15.054-000, S. J. Rio Preto, SP, Brasil }
\centerline{\email{c.pessoa@unesp.br} and \email{ronisio.ribeiro@unesp.br}}}

\medskip

\bigskip

\begin{quote}{\normalfont\fontsize{8}{10}\selectfont
{\bfseries Abstract.}
In this paper, we study the existence of limit cycles in continuous and discontinuous planar piecewise linear Hamiltonian differential system with two or three zones separated by straight lines and such that the linear systems that define the piecewise one have isolated singular points, i.e. centers or saddles. In this case, we show that if the planar piecewise linear Hamiltonian differential system is either continuous or discontinuous with two zones, then it has no limit cycles. Now, if the planar piecewise linear Hamiltonian differential system is discontinuous with three zones, then it has at most one limit cycle, and there are examples with one limit cycle. More precisely, without taking into account the position of the singular points in the zones, we present examples with the unique limit cycle for all possible combinations of saddles and centers.
\par}
\end{quote}

\section{Introduction}\label{sec:01}
The first works on piecewise differential systems appeared in the 1930s, see  \cite{And}. This class of systems is very important due to numerous applications, for example in control theory, mechanics, electrical circuits, neurobiology, etc (see for instance the book \cite{di}). Recently, this subject has piqued the attention of researchers in qualitative theory of differential equations and numerous studies about this topic have arisen in the literature (see \cite{Cas, Hua2, Li2, Li3, Mer}).
 
Piecewise linear differential systems are an interesting class of piecewise differential systems and, unlike the smooth case,  have a rich dynamic that is far from being fully understood. In addition to numerous applications in various areas of knowledge. In 1990 Lum and Chua \cite{Chu} conjectured that a continuous piecewise differential systems in the plane with two zones has at most one limit cycle. In 1998 this conjecture was proved by Freire, Ponce, Rodrigo and Torres in \cite{Fre}. The problem becomes more complicated when we have three zones. Conditions for non existence and existence of one, two or three limit cycles have been obtained, see \cite{Fre1, Lli4, Pon}. However, the maximum number of limit cycles, as far as we know, is not yet known.

In the discontinuous case, the maximum number of limit cycles is not known even in the simplest case, i.e. for piecewise linear differential systems with two zones separated by a straight line. However, important partial results about this problem have been obtained, see for instance \cite{Bra, Buz, Fre3, Fre2, Fre4, Hua, Li, Lli2, Wan2,  Wan1, Wan3}, among other papers. In summary, the results about the number of limit cycles of discontinuous piecewise linear differential systems with two zones separated by a straight line are given in Table \ref{tab1}. In that table the symbol “---” indicates that those cases appear repeated in the table.

\begin{table}[h!]\label{tab1}
\begin{tabular}{|c|c|c|c|c|c|}
\hline
       & Focus    & Center      & Saddle       & Node \\ \hline
Focus  & 3        & 2*           & 3           & 3    \\ \hline
Center & ---      & 0*           & 1*          & 1*   \\ \hline
Saddle & ---      & ---          & 2           & 2    \\ \hline
Node   & ---      & ---          & ---         & 2    \\ \hline
\end{tabular}
\caption{Lower bounds (Upper bounds*) of the maximum limit cycles of discontinuous piecewise linear differential systems with two zones separated by a straight line.}
\label{tab1}
\end{table}

If the curve between two linear zones is not a straight line is possible to obtain as many cycles as you want, see \cite{Bra2}. Exact number of limit cycles, for discontinuous piecewise linear systems with two zones separated by a straight line, were obtained in particular cases. Llibre and Teixeira \cite{Lli1} proved that if the linear systems, that define the piecewise one, has no singular point, then it has at most one limit cycle. Medrado and Torregrossa \cite{Med} proved that if the straight line has only crossing sewing points and the piecewise linear system has only a monodronic singular point on it, then the system has at most one limit cycle. 

There are few papers on discontinuous piecewise linear systems with three zones separated by two straight lines (see \cite{Mel,Lli3,Xio2,Xio}). In \cite{Lli3}, Llibre and Teixeira study the existence of limit cycles for continuous and discontinuous planar piecewise linear differential system with three zones separated by two parallel straight lines and such that the linear systems involved have a unique singular point which are centers. More precisely, in the continuous case, they prove that the piecewise system has no limit cycles. Now, in the discontinuous case, the piecewise system has at most one limit cycle and there are examples with one. Mello, Llibre and Fonseca, in \cite{Mel}, propose a mix of \cite{Lli1} and \cite{Lli3}. They proved that a piecewise linear Hamiltonian systems with three zones separated by two parallel straight lines without singular points have at most one crossing limit cycle.

In this paper, we contribute along these lines, that is, we are interested in studying the existence and the number of limit cycles of a piecewise linear differential systems with two or three zones in the plane with the following hypothesis:
\begin{itemize}
\item[{\rm (H1)}] The separations curves are straight lines.
\item[{\rm (H2)}] The vector fields which define the piecewise one are linear.
\item[{\rm (H3)}] The vector fields which define the piecewise one are Hamiltonian.
\item[{\rm (H4)}] The vector fields which define the piecewise one have isolated singularities. 
\end{itemize} 
Note that, hypothesis (H2), (H3) and (H4) imply that the singular points of the  linear systems that define the piecewise differential systems are saddles or centers.

We can classify the systems that satisfy the above hypotheses according to the configuration of their singular points. Thus, denoting the centers by the capital letter C and by S the saddles, in the case of two zones we have systems of the type CC, SC and SS. This is, CC indicates that the singular points of the linear systems that define the piecewise differential system are centers and so on.  Following this idea, for three zones, we have the following six class of piecewise linear Hamiltonian systems: CCC, SCC, SCS, CSC, SSS and SSC.

Assuming hypotheses (H1)--(H4), the main results in this paper are the follows:

\begin{theorem}\label{thm:01}
A continuous planar piecewise linear Hamiltonian differential system with two zones separated by a straight line and such that the linear systems that define it have isolated singular points, i.e. centers or saddles, has no limit cycles.
\end{theorem}

\begin{theorem}\label{thm:02}
A continuous planar piecewise linear Hamiltonian differential system with three zones separated by two parallel straight lines and such that the linear systems that define it have isolated singular points, i.e. centers or saddles, has no limit cycles.
\end{theorem}
 
\begin{theorem}\label{thm:03}
A discontinuous planar piecewise linear Hamiltonian differential system with two zones separated by a straight line and such that the linear systems that define it have isolated singular points, i.e. centers or saddles, has no limit cycles.
\end{theorem}

Theorem \ref{thm:03} has been proved in the literature and we have included it here just for the sake of completeness. See the proof of Theorem 1 from \cite{Lli2} for the case where one of the linear systems has a center and the other has a center or saddle, and see the proof of Theorem 3.4 from \cite{Hua3} for the case where the linear systems has isolated saddles.

\begin{theorem}\label{thm:04}
A discontinuous planar piecewise linear Hamiltonian differential system with three zones separated by two parallel straight lines and such that the linear systems that define it have isolated singular points, i.e. centers or saddles, has at most a limit cycles.
\end{theorem}

Theorems \ref{thm:01}, \ref{thm:02} and \ref{thm:04} have been proved for the particular case in which the linear systems that define the piecewise one has only isolated centers, see \cite{Lli3}. Theorem \ref{thm:02} has also been proved for the particular case SCS, see the proof of Lemma 11 from \cite{Pu}. For the other possibilities, as far as we know, the results of Theorems 1, 2 and 4 are new.  

\medskip

The paper is organized as follows. In section \ref{Sec:02} we introduce the basic definitions and results. In Section \ref{Sec:03} we prove Theorems \ref{thm:01}-\ref{thm:04}. Examples of discontinuous planar piecewise linear Hamiltonian differential system with three zones separated by two parallel straight lines such that the linear systems that define it have isolated singular points are analyzed in section \ref{Sec:04}. That is, we give examples of piecewise linear Hamiltonian systems of type CCC, SCC, SCS, CSC, SSS and SSC with exactly one limit cycle.

\section{Preliminary results}\label{Sec:02}
 
In this section, we will present the basic concepts that we need to prove the main results of this paper.

Let $h_i:\R^2\rightarrow\R$,  $i=C,L,R$, be the function  $h_{\scriptscriptstyle C}(x,y)=x$, $h_{\scriptscriptstyle L}(x,y)=x+1$ and $h_{\scriptscriptstyle R}(x,y)=x-1$. The {\it switching curve} $\Sigma_{\scriptscriptstyle C}$ of a piecewise linear system with two zones in the plane is defined as
$$\Sigma_{\scriptscriptstyle C}=h_{\scriptscriptstyle C}^{-1}(0)=\{(x,y)\in\R^2:x=0\}.$$
This straight line decomposes the plane in two regions 
$$R_{\scriptscriptstyle L}=\{(x,y)\in\R^2:x<0\}\quad\text{and}\quad R_{\scriptscriptstyle R}=\{(x,y)\in\R^2:x>0\}.$$

Assuming the hypothesis (H2) and (H3), the piecewise linear Hamiltonian vector field with two zones is given by
\begin{equation}\label{sys1}
   \begin{aligned}
      X_{\scriptscriptstyle L}(x,y)&=(a_{\scriptscriptstyle L}x+b_{\scriptscriptstyle L}y+\alpha_{\scriptscriptstyle L},c_{\scriptscriptstyle L}x-a_{\scriptscriptstyle L}y+\beta_{\scriptscriptstyle L}),\quad x\leq 0, \\
			X_{\scriptscriptstyle R}(x,y)&=(a_{\scriptscriptstyle R}x+b_{\scriptscriptstyle R}y+\alpha_{\scriptscriptstyle R},c_{\scriptscriptstyle R}x-a_{\scriptscriptstyle R}y+\beta_{\scriptscriptstyle R}),\quad x>0.
  \end{aligned}
\end{equation}
Note that the Hamiltonians functions that determine the vector field (\ref{sys1}) are
\begin{equation}\label{ham1}  
\begin{aligned}
      H_{\scriptscriptstyle L}(x,y)&=\frac{b_{\scriptscriptstyle L}}{2}y^2-\frac{c_{\scriptscriptstyle L}}{2}x^2+a_{\scriptscriptstyle L}xy+\alpha_{\scriptscriptstyle L}y-\beta_{\scriptscriptstyle L}x,\quad x\leq 0, \\
			H_{\scriptscriptstyle R}(x,y)&=\frac{b_{\scriptscriptstyle R}}{2}y^2-\frac{c_{\scriptscriptstyle R}}{2}x^2+a_{\scriptscriptstyle R}xy+\alpha_{\scriptscriptstyle R}y-\beta_{\scriptscriptstyle R}x,\quad x>0.
 \end{aligned}
\end{equation}  
 
For the case with three zones, the switching curves $\Sigma_{\scriptscriptstyle L}$ and $\Sigma_{\scriptscriptstyle R}$ are given by
$$\Sigma_{\scriptscriptstyle L}=h_{\scriptscriptstyle L}^{-1}(0)=\{(x,y)\in\R^2:x=-1\},$$
and
$$\Sigma_{\scriptscriptstyle R}=h_{\scriptscriptstyle R}^{-1}(0)=\{(x,y)\in\R^2:x=1\}.$$
This straight line decomposes the plane in three regions 
$$R_{\scriptscriptstyle L}=\{(x,y)\in\R^2:x<-1\},\quad R_{{\scriptscriptstyle C}}=\{(x,y)\in\R^2:-1<x<1\},$$ 
and 
$$R_{{\scriptscriptstyle R}}=\{(x,y)\in\R^2:x>1\}.$$
Assuming the hypothesis (H2) and (H3), the piecewise linear Hamiltonian vector field with three zones is give by
\begin{equation}\label{sys2}
   \begin{aligned}
      X_{\scriptscriptstyle L}(x,y)&=(a_{\scriptscriptstyle L}x+b_{\scriptscriptstyle L}y+\alpha_{\scriptscriptstyle L},c_{\scriptscriptstyle L}x-a_{\scriptscriptstyle L}y+\beta_{\scriptscriptstyle L}),\quad x \leq -1, \\
      X_{\scriptscriptstyle C}(x,y)&=(a_{\scriptscriptstyle C}x+b_{\scriptscriptstyle C}y+\alpha_{\scriptscriptstyle C},c_{\scriptscriptstyle C}x-a_{\scriptscriptstyle C}y+\beta_{\scriptscriptstyle C}),\quad -1\leq x\leq 1, \\
			X_{\scriptscriptstyle R}(x,y)&=(a_{\scriptscriptstyle R}x+b_{\scriptscriptstyle R}y+\alpha_{\scriptscriptstyle R},c_{\scriptscriptstyle R}x-a_{\scriptscriptstyle R}y+\beta_{\scriptscriptstyle R}),\quad x \geq 1.
  \end{aligned}
\end{equation}
The Hamiltonians functions that determine the vector field (\ref{sys2}) are
\begin{equation}\label{ham2}  
\begin{aligned}
      H_{\scriptscriptstyle L}(x,y)&=\frac{b_{\scriptscriptstyle L}}{2}y^2-\frac{c_{\scriptscriptstyle L}}{2}x^2+a_{\scriptscriptstyle L}xy+\alpha_{\scriptscriptstyle L}y-\beta_{\scriptscriptstyle L}x,\quad x\leq -1, \\
      H_{\scriptscriptstyle C}(x,y)&=\frac{b_{\scriptscriptstyle C}}{2}y^2-\frac{c_{\scriptscriptstyle C}}{2}x^2+a_{\scriptscriptstyle C}xy+\alpha_{\scriptscriptstyle C}y-\beta_{\scriptscriptstyle C}x,\quad -1\leq x\leq 1, \\
			H_{\scriptscriptstyle R}(x,y)&=\frac{b_{\scriptscriptstyle R}}{2}y^2-\frac{c_{\scriptscriptstyle R}}{2}x^2+a_{\scriptscriptstyle R}xy+\alpha_{\scriptscriptstyle R}y-\beta_{\scriptscriptstyle R}x,\quad x\geq 1.
 \end{aligned}
\end{equation}

We will use the vector field $X_{\scriptscriptstyle L}$ and the switching curve $\Sigma_{\scriptscriptstyle L}$ in the next definitions. However, we can easily adapt the definitions to the vector fields $X_{\scriptscriptstyle C}$ and $X_{\scriptscriptstyle R}$ and the switching curves $\Sigma_{\scriptscriptstyle C}$ and $\Sigma_{\scriptscriptstyle R}$.

The derivative of function $h_{\scriptscriptstyle L}$ in the direction of the vector field $X_{\scriptscriptstyle L}$, i.e., the expression
$$X_{\scriptscriptstyle L} h_{\scriptscriptstyle L}(p)=\langle  X_{\scriptscriptstyle L}(p),\nabla h_{\scriptscriptstyle L}(p)\rangle,$$ 
where $\left\langle \cdot,\cdot\right\rangle$ is the usual inner product in $\R^2$,  characterize the contact between the vector field $X_{\scriptscriptstyle L}$ and the switching curve $\Sigma_{\scriptscriptstyle L}$. 

We distinguish the followings subsets of $\Sigma_{\scriptscriptstyle L}$ (the same for $\Sigma_{\scriptscriptstyle C}$ and $\Sigma_{\scriptscriptstyle R}$)
\begin{itemize}
\item Crossing set:
$$\Sigma_{\scriptscriptstyle L}^{c}=\{p\in \Sigma_{\scriptscriptstyle L}:X_{\scriptscriptstyle L} h_{\scriptscriptstyle L}(p) \cdot X_{\scriptscriptstyle C} h_{\scriptscriptstyle L}(p)>0\};$$
\item Sliding set:
$$\Sigma_{\scriptscriptstyle L}^{s}=\{p\in \Sigma_{\scriptscriptstyle L}:X_{\scriptscriptstyle L} h_{\scriptscriptstyle L}(p)>0, X_{\scriptscriptstyle C} h_{\scriptscriptstyle L}(p)<0\};$$
\item Escaping set:
$$\Sigma_{\scriptscriptstyle L}^{e}=\{p\in \Sigma_{\scriptscriptstyle L}:X_{\scriptscriptstyle L} h_{\scriptscriptstyle L}(p)<0, X_{\scriptscriptstyle C} h_{\scriptscriptstyle L}(p)>0\}.$$
\end{itemize}
In a piecewise vector field with two or three zones in the plane, the limit cycles can be of two types: sliding limit cycles or crossing limit cycles; the first one contain some segment of sliding or escaping sets, and the second one does not contain any segments of sliding or escaping sets. In this paper, we only study the crossing limit cycles. In what follows, when we talk about limit cycles, we are talking about crossing limit cycles.

A piecewise vector field with two zones in the plane is called continuous if 
\begin{equation*}
X_{\scriptscriptstyle L}(p)= X_{\scriptscriptstyle R}(p),\quad \forall p \in \Sigma_{\scriptscriptstyle C}.
\end{equation*}
Otherwise, it is called discontinuous. Similarly, a piecewise vector field with three zones in the plane is called continuous if 
\begin{equation*}  
\begin{aligned}
      X_{\scriptscriptstyle L}(p)&= X_{\scriptscriptstyle C}(p),\quad \forall p \in \Sigma_{\scriptscriptstyle L} \quad \text{e} \\
      X_{\scriptscriptstyle C}(q)&= X_{\scriptscriptstyle R}(q),\quad \forall q \in \Sigma_{\scriptscriptstyle R}.
 \end{aligned}
\end{equation*}
Otherwise, it is called discontinuous.

\medskip

\section{Proof of Theorems \ref{thm:01}-\ref{thm:04}}\label{Sec:03}
This section is devoted to present the proof of main results.

\medskip

\noindent \textit{Proof of Theorem \ref{thm:01}}. Consider a continuous piecewise linear Hamiltonian vector field with two zones separated by a straight line, such that the linear vector fields, that define it, have isolated singular points. That is, we have piecewise vector field (\ref{sys1}), with $a_{\scriptscriptstyle i}^2+b_{\scriptscriptstyle i}c_{\scriptscriptstyle i}\ne 0$, for $i=L,R$, and due to continuity
$$ X_{\scriptscriptstyle R}(0,y)= X_{\scriptscriptstyle L}(0,y),\quad\forall y\in\R.$$
This equality implies that
$$a_{\scriptscriptstyle R}=a_{\scriptscriptstyle L}=a,\quad b_{\scriptscriptstyle R}=b_{\scriptscriptstyle L}=b,\quad \alpha_{\scriptscriptstyle R}=\alpha_{\scriptscriptstyle L}=\alpha\quad\text{and}\quad\beta_{\scriptscriptstyle R}=\beta_{\scriptscriptstyle L}=\beta.$$ 

If the piecewise linear vector field has a periodic orbit, then it intersects the straight line $x=0$ at two points, $(0,y_0)$ and $(0,y_1)$, with $y_1<y_0$, satisfying 
\begin{equation*}  
\begin{aligned}
      H_{\scriptscriptstyle R}(0,y_1)&= H_{\scriptscriptstyle R}(0,y_0), \\
      H_{\scriptscriptstyle L}(0,y_0)&= H_{\scriptscriptstyle L}(0,y_1),
 \end{aligned}
\end{equation*}
where  $H_{\scriptscriptstyle L}$ and $H_{\scriptscriptstyle R}$ are given by (\ref{ham1}). More precisely, we have the equations
\begin{equation*}  
\begin{aligned}
      -\frac{1}{2} (y_0 - y_1) (b (y_0 + y_1) + 2 \alpha)&=0, \\
      \frac{1}{2} (y_0 - y_1) (b (y_0 + y_1) + 2 \alpha)&=0.
 \end{aligned}
\end{equation*}
As $y_1<y_0$, if $b=0$ and $\alpha\ne 0$, then the system above has no solution. If $b=\alpha=0$, then the system has infinity solutions.
If $b\ne 0$ the solution $(y_0,y_1)$ of the above system with $y_1<y_0$ satisfies $y_0=-(by_1+2\alpha)/b$, with arbitrary $y_1$. Therefore, the piecewise linear vector field (\ref{sys1}) has no periodic orbits or has a continuum of periodic orbits, and consequently, it has no limit cycle.
$\hfill\square$

\bigskip

\noindent \textit{Proof of Theorem \ref{thm:02}}. Consider a continuous piecewise linear Hamiltonian vector field with three zones separated by two parallel straight lines, such that the linear vector fields, that define it, have isolated singular points. That is, we have  piecewise vector fields (\ref{sys2}), with $a_{\scriptscriptstyle i}^2+b_{\scriptscriptstyle i}c_{\scriptscriptstyle i}\ne 0$, for $i=L,C,R$, and due to continuity
$$ X_{\scriptscriptstyle R}(1,y)= X_{\scriptscriptstyle C}(1,y)\quad\text{and}\quad X_{\scriptscriptstyle C}(-1,y)= X_{\scriptscriptstyle L}(-1,y),\quad\forall y\in\R.$$
These equalities imply that
$$a_{\scriptscriptstyle R}=a_{\scriptscriptstyle C}=a_{\scriptscriptstyle L}=a,\quad b_{\scriptscriptstyle R}=b_{\scriptscriptstyle C}=b_{\scriptscriptstyle L}=b,\quad\alpha_{\scriptscriptstyle R}=\alpha_{\scriptscriptstyle C}=\alpha_{\scriptscriptstyle L}=\alpha$$
and
$$\beta_{\scriptscriptstyle R}-\beta_{\scriptscriptstyle C}-c_{\scriptscriptstyle C}+c_{\scriptscriptstyle R}=\beta_{\scriptscriptstyle L}-\beta_{\scriptscriptstyle C}-c_{\scriptscriptstyle L}+c_{\scriptscriptstyle C}=0.$$
If the piecewise linear vector field has a periodic orbit, then it intersects the straight lines $x=\pm 1$ at four points, $(1,y_0)$, $(1,y_1)$, with $y_1<y_0$, and $(-1,y_2)$, $(-1,y_3)$, with $y_2<y_3$, respectively, satisfying
\begin{equation}\label{eqc5}  
\begin{aligned}
      &H_{\scriptscriptstyle R}(1,y_1)= H_{\scriptscriptstyle R}(1,y_0), \\
			&H_{\scriptscriptstyle C}(1,y_0)= H_{\scriptscriptstyle C}(-1,y_3), \\
			&H_{\scriptscriptstyle L}(-1,y_3)= H_{\scriptscriptstyle L}(-1,y_2), \\
      &H_{\scriptscriptstyle C}(-1,y_2)= H_{\scriptscriptstyle C}(1,y_1),
 \end{aligned}
\end{equation}
where  $H_{\scriptscriptstyle L}$, $H_{\scriptscriptstyle C}$ and $H_{\scriptscriptstyle R}$ are given by (\ref{ham2}). More precisely, we have the equations 
\begin{eqnarray}  
     &&-\frac{1}{2} (y_0 - y_1) (b (y_0 + y_1) + 2 (a + \alpha))=0, \label{c1}\\
     &&a (y_0 + y_3) + \frac{1}{2} (y_0 - y_3) (b (y_0 + y_3) + 2 \alpha) - 2 \beta_{\scriptscriptstyle C}=0,\label{c2} \\
		 &&-\frac{1}{2} (y_2 - y_3) (b (y_2 + y_3) -2 (a - \alpha))=0, \label{c3}\\
		 &&-a (y_1 + y_2) - \frac{1}{2} (y_1 - y_2) (b (y_1 + y_2) + 2 \alpha) + 2 \beta_{\scriptscriptstyle C}=0 \label{c4}.
\end{eqnarray}
As $y_1<y_0$, $y_2<y_3$ and $a^2+bc_{\scriptscriptstyle i}\ne 0$, for $i=L,C,R$, if either $b=0$ and $a+\alpha\neq 0$ or $b=0$ and $a+\alpha= 0$  the above system has no solutions. If $b\ne 0$, as $y_1<y_0$ and $y_2<y_3$, from equation (\ref{c1}) we can obtain $y_0$ as a function of $y_1$, i.e.
\begin{equation}\label{eqc6}
y_0 =\frac{- b y_1 - 2(a + \alpha)}{b}.
\end{equation}
Now, from equation (\ref{c3}) we can obtain $y_2$ as a function of $y_3$, i.e.
\begin{equation}\label{eqc7}
y_2 =\frac{- b y_3 - 2(\alpha-a)}{b}.
\end{equation}
Substituting (\ref{eqc6}) and (\ref{eqc7}) in equations (\ref{c2}) and (\ref{c4}), respectively, we obtain a solution $(y_0,y_1,y_2,y_3)$ of the system (\ref{eqc5}) satisfying $y_1<y_0$ and $y_2<y_3$, given by $(\varphi_1(y_1),y_1,\varphi_2(y_1),\varphi_3(y_1))$, where
\begin{equation*}
\begin{aligned}
      &\varphi_1(y_1)=\frac{-b y_1-2(a +\alpha)}{b}, \\
	  &\varphi_2(y_1)= \frac{a-\alpha+\sqrt{a^2+2 a (b y_1-\alpha)+(b y_1+\alpha)^2-4 b \beta_{\scriptscriptstyle C}}}{b}, \\
	  &\varphi_3(y_1)= \frac{a-\alpha-\sqrt{a^2+2 a (b y_1-\alpha)+(b y_1+\alpha)^2-4 b \beta_{\scriptscriptstyle C}}}{b},
 \end{aligned}
\end{equation*}
with arbitrary $y_1$. Therefore, the piecewise linear vector field (\ref{sys2}) has no periodic orbits or has a continuum of periodic orbits, and consequently, it has no limit cycle.
$\hfill\square$

\bigskip

\noindent\textit{Proof of Theorem \ref{thm:03}}. Consider a discontinuous piecewise linear Hamiltonian vector field with two zones separated by a straight line, such that the linear vector fields, that define it, have isolated singular points. That is, we have piecewise vector field (\ref{sys1}), with $a_{\scriptscriptstyle i}^2+b_{\scriptscriptstyle i}c_{\scriptscriptstyle i}\ne 0$, for $i=L,R$. If the piecewise linear vector field has a periodic orbit, then it intersects the straight line $x=0$ at two points, $(0,y_0)$ and $(0,y_1)$, with $y_1<y_0$, satisfying
\begin{equation*}  
\begin{aligned}
      H_{\scriptscriptstyle R}(0,y_1)&= H_{\scriptscriptstyle R}(0,y_0), \\
      H_{\scriptscriptstyle L}(0,y_0)&= H_{\scriptscriptstyle L}(0,y_1),
 \end{aligned}
\end{equation*}
where  $H_{\scriptscriptstyle L}$ and $H_{\scriptscriptstyle R}$ are given by (\ref{ham1}). More precisely, we have the equations
\begin{equation*}  
\begin{aligned}
      -\frac{1}{2} (y_0 - y_1) (b_{\scriptscriptstyle R} (y_0 + y_1) + 2 \alpha_{\scriptscriptstyle R})&=0, \\
      \frac{1}{2} (y_0 - y_1) (b_{\scriptscriptstyle L} (y_0 + y_1) + 2 \alpha_{\scriptscriptstyle L})&=0.
 \end{aligned}
\end{equation*}
As $y_1<y_0$, if $b_{\scriptscriptstyle R}=0$ and $\alpha_{\scriptscriptstyle R}\ne 0$ or $b_{\scriptscriptstyle L}=0$ and $\alpha_{\scriptscriptstyle L}\ne 0$  the above system has no solutions. If $b_{\scriptscriptstyle R}=\alpha_{\scriptscriptstyle R}=0$ and $b_{\scriptscriptstyle L}\ne0$ the solution $(y_0,y_1)$ of the above system with $y_1<y_0$ satisfies $y_0=-(b_{\scriptscriptstyle L}y_1+2\alpha_{\scriptscriptstyle L})/b_{\scriptscriptstyle L}$, with arbitrary $y_1$. 
If $b_{\scriptscriptstyle L}=\alpha_{\scriptscriptstyle L}=0$ and $b_{\scriptscriptstyle R}\ne0$ the solution $(y_0,y_1)$ of the above system with $y_1<y_0$ satisfies $y_0=-(b_{\scriptscriptstyle R}y_1+2\alpha_{\scriptscriptstyle R})/b_{\scriptscriptstyle R}$, with arbitrary $y_1$. If $b_{\scriptscriptstyle L}b_{\scriptscriptstyle R}\ne 0$, then the above system  has a solution $(y_0,y_1)$ with $y_1<y_0$ only when $b_{\scriptscriptstyle L}=b_{\scriptscriptstyle R}=b$ and $\alpha_{\scriptscriptstyle L}=\alpha_{\scriptscriptstyle R}=\alpha$. Moreover, $y_0=(-by_1+2\alpha)/b$ with arbitrary $y_1$.
If $b_{\scriptscriptstyle R}=b_{\scriptscriptstyle L}=\alpha_{\scriptscriptstyle R}=\alpha_{\scriptscriptstyle L}=0$, then the system has infinity solutions.
Therefore, the piecewise linear vector field (\ref{sys1}) has no periodic orbits or has a continuum of periodic orbits, and consequently, it has no limit cycle.
$\hfill\square$

\bigskip

\noindent \textit{Proof of Theorem \ref{thm:04}}. Consider a discontinuous piecewise linear Hamiltonian vector field with three zones separated by two parallel  straight lines, such that the linear vector fields, that define it, have isolated singular points. That is, we have  piecewise vector fields (\ref{sys2}), with $-a_{\scriptscriptstyle i}^2-b_{\scriptscriptstyle i}c_{\scriptscriptstyle i}\ne 0$, for $i=L,C,R$. If the piecewise linear vector field has a periodic orbit, then it intersects the straight lines $x=\pm 1$ at four points, $(1,y_0)$, $(1,y_1)$, with $y_1<y_0$, and $(-1,y_2)$, $(-1,y_3)$, with $y_2<y_3$, respectively, satisfying
\begin{equation}\label{eq5}  
\begin{aligned}
      &H_{\scriptscriptstyle R}(1,y_1)= H_{\scriptscriptstyle R}(1,y_0), \\
			&H_{\scriptscriptstyle C}(1,y_0)= H_{\scriptscriptstyle C}(-1,y_3), \\
			&H_{\scriptscriptstyle L}(-1,y_3)= H_{\scriptscriptstyle L}(-1,y_2), \\
      &H_{\scriptscriptstyle C}(-1,y_2)= H_{\scriptscriptstyle C}(1,y_1),
 \end{aligned}
\end{equation}
where  $H_{\scriptscriptstyle L}$, $H_{\scriptscriptstyle C}$ and $H_{\scriptscriptstyle R}$ are given by (\ref{ham2}). More precisely, we have the equations 
\begin{eqnarray}  
     &&\frac{1}{2} (y_1 - y_0) (b_{\scriptscriptstyle R} (y_0 + y_1) + 2 (a_{\scriptscriptstyle R} + \alpha_{\scriptscriptstyle R}))=0, \label{s1}\\
     && \frac{1}{2} (y_0 - y_3) (b_{\scriptscriptstyle C} (y_0 + y_3) + 2 \alpha_{\scriptscriptstyle C}) - 2 \beta_{\scriptscriptstyle C}+a_{\scriptscriptstyle C} (y_0 + y_3)=0,\label{s2} \\
		 &&\frac{1}{2} (y_3 - y_2) (b_{\scriptscriptstyle L} (y_2 + y_3) - 2 (a_{\scriptscriptstyle L} - \alpha_{\scriptscriptstyle L}))=0, \label{s3}\\
		 &&\frac{1}{2} (y_2 - y_1) (b_{\scriptscriptstyle C} (y_1 + y_2) + 2 \alpha_{\scriptscriptstyle C}) + 2 \beta_{\scriptscriptstyle C}-a_{\scriptscriptstyle C} (y_1 + y_2)=0 \label{s4}.
\end{eqnarray} 
As $y_1<y_0$, $y_2<y_3$ and $a_{\scriptscriptstyle i}^2+b_{\scriptscriptstyle i}c_{\scriptscriptstyle i}\ne 0$, for $i=L,C,R$, if $b_{\scriptscriptstyle R}=0$ and $a_{\scriptscriptstyle R}+\alpha_{\scriptscriptstyle R}\ne 0$ or $b_{\scriptscriptstyle L}=0$ and $a_{\scriptscriptstyle L}-\alpha_{\scriptscriptstyle L}\ne 0$ or $b_{\scriptscriptstyle R}=a_{\scriptscriptstyle R}+\alpha_{\scriptscriptstyle R}=b_{\scriptscriptstyle L}=a_{\scriptscriptstyle L}-\alpha_{\scriptscriptstyle L}=b_{\scriptscriptstyle C}=\alpha_{\scriptscriptstyle C}-a_{\scriptscriptstyle C}=0$ or $b_{\scriptscriptstyle R}=a_{\scriptscriptstyle R}+\alpha_{\scriptscriptstyle R}=b_{\scriptscriptstyle C}=\alpha_{\scriptscriptstyle C}-a_{\scriptscriptstyle C}=0$ and $b_{\scriptscriptstyle L}\ne0$ or $b_{\scriptscriptstyle R}\ne 0$ and $b_{\scriptscriptstyle L}=a_{\scriptscriptstyle L}-\alpha_{\scriptscriptstyle L}=b_{\scriptscriptstyle C}=\alpha_{\scriptscriptstyle C}-a_{\scriptscriptstyle C}=0$ or $b_{\scriptscriptstyle R}b_{\scriptscriptstyle L}\ne 0$, $b_{\scriptscriptstyle C}=0$ and $b_{\scriptscriptstyle R} \alpha_{\scriptscriptstyle C} (a_{\scriptscriptstyle L} - \alpha_{\scriptscriptstyle L}) + a_{\scriptscriptstyle C} b_{\scriptscriptstyle R} (\alpha_{\scriptscriptstyle L} - a_{\scriptscriptstyle L}) + b_{\scriptscriptstyle L} (a_{\scriptscriptstyle R} + \alpha_{\scriptscriptstyle R})(a_{\scriptscriptstyle C} + \alpha_{\scriptscriptstyle C}) + 2 b_{\scriptscriptstyle L} b_{\scriptscriptstyle R} \beta_{\scriptscriptstyle C}\ne 0$ the above system has no solutions.
If $b_{\scriptscriptstyle R}=a_{\scriptscriptstyle R}+\alpha_{\scriptscriptstyle R}=b_{\scriptscriptstyle L}=a_{\scriptscriptstyle L}-\alpha_{\scriptscriptstyle L}=b_{\scriptscriptstyle C}=0$ and $\alpha_{\scriptscriptstyle C}-a_{\scriptscriptstyle C}\ne 0$ or $b_{\scriptscriptstyle R}=a_{\scriptscriptstyle R}+\alpha_{\scriptscriptstyle R}=b_{\scriptscriptstyle L}=a_{\scriptscriptstyle L}-\alpha_{\scriptscriptstyle L}=0$ and $b_{\scriptscriptstyle C}\ne0$ or $b_{\scriptscriptstyle R}=a_{\scriptscriptstyle R}+\alpha_{\scriptscriptstyle R}=b_{\scriptscriptstyle C}=0$, $\alpha_{\scriptscriptstyle C}-a_{\scriptscriptstyle C}\ne0$ and $b_{\scriptscriptstyle L}\ne0$ or $b_{\scriptscriptstyle R}=a_{\scriptscriptstyle R}+\alpha_{\scriptscriptstyle R}=0$ and $b_{\scriptscriptstyle L}b_{\scriptscriptstyle C}\ne0$ or $b_{\scriptscriptstyle R}\ne0$, $b_{\scriptscriptstyle L}=a_{\scriptscriptstyle L}-\alpha_{\scriptscriptstyle L}=b_{\scriptscriptstyle C}=0$ and $\alpha_{\scriptscriptstyle C}-a_{\scriptscriptstyle C}\ne0$ or $b_{\scriptscriptstyle R}b_{\scriptscriptstyle C}\ne0$ and $b_{\scriptscriptstyle L}=a_{\scriptscriptstyle L}-\alpha_{\scriptscriptstyle L}=0$ or $b_{\scriptscriptstyle R}b_{\scriptscriptstyle L}\ne0$, $b_{\scriptscriptstyle C}=0$ and  $b_{\scriptscriptstyle R} \alpha_{\scriptscriptstyle C} (a_{\scriptscriptstyle L} - \alpha_{\scriptscriptstyle L}) + a_{\scriptscriptstyle C} b_{\scriptscriptstyle R} (\alpha_{\scriptscriptstyle L} - a_{\scriptscriptstyle L}) + b_{\scriptscriptstyle L} (a_{\scriptscriptstyle R} + \alpha_{\scriptscriptstyle R})(a_{\scriptscriptstyle C} + \alpha_{\scriptscriptstyle C}) + 2 b_{\scriptscriptstyle L} b_{\scriptscriptstyle R} \beta_{\scriptscriptstyle C}= 0$ then the system (\ref{s1})--(\ref{s4}) has infinity solutions.  
If $b_{\scriptscriptstyle L}b_{\scriptscriptstyle C}b_{\scriptscriptstyle R}\ne 0$, from equation (\ref{s1}), we can obtain $y_0$ as a function of $y_1$, i.e.
\begin{equation}\label{eq6}
y_0 =\frac{- b_{\scriptscriptstyle R} y_1 - 2 (a_{\scriptscriptstyle R} + \alpha_{\scriptscriptstyle R})}{b_{\scriptscriptstyle R}}.
\end{equation}
Now, from equation (\ref{s3}), we can obtain $y_2$ as a function of $y_3$, i.e.
\begin{equation}\label{eq7}
y_2 =\frac{- b_{\scriptscriptstyle L} y_3 - 2 (\alpha_{\scriptscriptstyle L} - a_{\scriptscriptstyle L})}{b_{\scriptscriptstyle L}}.
\end{equation}
Substituting (\ref{eq6}) and (\ref{eq7}) in equations (\ref{s2}) and (\ref{s4}), respectively, we obtain the equations of two hyperbolas in the $y_1y_3$ plane, given by
\begin{equation}\label{eq8}
\begin{aligned}
     &\frac{(y_1-A)^2}{K}-\frac{(y_3-B)^2}{K}-C=0, \\
     &\frac{(y_1-D)^2}{K}-\frac{(y_3-E)^2}{K}-C=0,
\end{aligned}
\end{equation}
with

$$K=\frac{2}{b_{\scriptscriptstyle C}},\quad A=\frac{b_{\scriptscriptstyle R}(a_{\scriptscriptstyle C}+\alpha_{\scriptscriptstyle C})-2b_{\scriptscriptstyle C}(a_{\scriptscriptstyle R}+\alpha_{\scriptscriptstyle R})}{b_{\scriptscriptstyle C} b_{\scriptscriptstyle R}},$$\vspace{0.1cm}
$$B=\frac{a_{\scriptscriptstyle C}-\alpha_{\scriptscriptstyle C}}{b_{\scriptscriptstyle C}}, \quad C=\frac{2(a_{\scriptscriptstyle C} \alpha_{\scriptscriptstyle C}+b_{\scriptscriptstyle C} \beta_{\scriptscriptstyle C})}{b_{\scriptscriptstyle C}},$$\vspace{0.1cm}
$$D=-\frac{(a_{\scriptscriptstyle C}+\alpha_{\scriptscriptstyle C})}{b_{\scriptscriptstyle C}}\quad\text{and}\quad E=\frac{b_{\scriptscriptstyle L}(\alpha_{\scriptscriptstyle C}-a_{\scriptscriptstyle C})-2b_{\scriptscriptstyle C}(\alpha_{\scriptscriptstyle L}-a_{\scriptscriptstyle L})}{b_{\scriptscriptstyle C} b_{\scriptscriptstyle L}}.$$
Note that the system (\ref{eq8}) is equivalent to the system
\begin{equation}\label{eq10} 
	\begin{aligned}
     &y_1^2-2Ay_1+A^2-y_3^2+2By_3-B^2-KC=0, \\
     &2(A-D)y_1+2(E-B)y_3+D^2-E^2+B^2-A^2=0.
    \end{aligned} 
\end{equation}
The system above eventually could have infinity solutions $(y_1,y_3)$, for instance when $A=D$ and $B=E$. In this case, the piecewise linear vector field (\ref{sys2}) has a continuum of periodic orbits, and consequently, it has no limit cycle. Suppose that system (\ref{eq10}) has finitely many solutions.
According to Bezout's Theorem, if a system of polynomial equations has
finitely many solutions, then the number of its solutions is at most the product of the degrees of the polynomials, that for system (\ref{eq10}) is two. Therefore, the two hyperbolas above intersect at most two points. Note that, by (\ref{s1})--(\ref{s4}), if $(y_0,y_1,y_2,y_3)$ is solution of the system (\ref{eq5}) then $(y_1,y_0,y_3,y_2)$ is also a solution. However, for $y_1<y_0$ and $y_2<y_3$ we have at most a single solution. Therefore, the piecewise linear vector field (\ref{sys2}) can have at most one limit cycle.
$\hfill\square$

\medskip

\section{Exemples}\label{Sec:04}
In this section, we will give some examples of discontinuous planar piecewise linear Hamiltonian differential system with three zones separated by two parallel straight lines with one limit cycle, such that the linear systems that define it have isolated singular points. That is, we given examples of piecewise linear Hamiltonian systems of type CCC, SCC, SCS, CSC, SSS and SSC with exactly one limit cycle. In \cite{Lli3}, the authors presented an example of a discontinuous piecewise linear differential system of type CCC with exactly one limit cycle. Here we will show another example for this case.

\begin{example} (Case CCC)
\label{CCC}
Consider the discontinuous planar piecewise linear Hamiltonian vector field \eqref{sys2} with $a_{\scriptscriptstyle L}=4$, $b_{\scriptscriptstyle L}=8$, $\alpha_{\scriptscriptstyle L}=3/2$, $c_{\scriptscriptstyle L}=-5/2$, $\beta_{\scriptscriptstyle L}=11/4$, $a_{\scriptscriptstyle C}=0$, $b_{\scriptscriptstyle C}=2$, $\alpha_{\scriptscriptstyle C}=\beta_{\scriptscriptstyle C}=2/3$, $c_{\scriptscriptstyle C}=-2$, $a_{\scriptscriptstyle R}=4$, $b_{\scriptscriptstyle R}=2$, $c_{\scriptscriptstyle R}=-10$ and $\alpha_{\scriptscriptstyle R}=\beta_{\scriptscriptstyle R}=-4$.
%\begin{equation}\label{ex1}
%   \begin{aligned}
%      X_{\scriptscriptstyle L}(x,y)&=\Big(4x+8y+\frac{3}{2},-\frac{5x}{2} - 4 y + \frac{11}{4}\Big),\quad x \leq -1, \\
%      X_{\scriptscriptstyle C}(x,y)&=\Big(2 y + \frac{2}{3},- 2 x + \frac{2}{3}\Big),\quad -1\leq x\leq 1, \\
%	  X_{\scriptscriptstyle R}(x,y)&=( 4 x + 2 y - 4, - 10 x - 4 y - 4),\quad x \geq 1.
%  \end{aligned}
%\end{equation}
The eigenvalues of the linear part of $X_i$, $i=L, C, R$,  from \eqref{sys2} for this case, are $\pm 2i$, $\pm 2i$ and $\pm 2i$, respectively, i.e. we have three centers. 
%The Hamiltonian functions that determine the systems (\ref{ex1}) are
%\begin{equation*}  
%\begin{aligned}
%      H_{\scriptscriptstyle L}(x,y)&=\frac{5x^2}{4} + 4y^2 + 4xy + \frac{3y}{2} - \frac{11x}{4},\quad x\leq -1, \\
%      H_{\scriptscriptstyle C}(x,y)&=x^2+y^2-\frac{2x}{3} + \frac{2y}{3},\quad -1\leq x\leq 1, \\
%	  H_{\scriptscriptstyle R}(x,y)&=5x^2 + y^2 + 4xy + 4x - 4y,\quad x\geq 1,
% \end{aligned}
%\end{equation*}
%respectively. 
Therefore, a candidate to limit cycle of vector field  \eqref{sys2}, in this case, correspond to the solution of system \eqref{eq5}, i.e.
\begin{equation*}\label{ex1s}  
\begin{aligned} 
     &(y_1-y_0)(y_1+y_0)=0,\\
     &\frac{1}{3}(y_0 (2 + 3 y_0)- y_3 (2 + 3 y_3) - 4)=0,\\
		 &\frac{1}{2}(y_3 - y_2) (8 (y_2 + y_3) - 5)=0,\\
		 &\frac{1}{3}(4 - y_1 (2 + 3 y_1) + y_2 (2 + 3 y_2))=0.
\end{aligned}
\end{equation*}
After some computations, the unique solution $(y_0,y_1,y_2,y_3)$ of the above system, satisfying the condition $y_1<y_0$ and $y_2<y_3$, is given by
$$\Bigg(\frac{31}{48}\sqrt{\frac{1259}{235}},-\frac{31}{48}\sqrt{\frac{1259}{235}},\frac{5}{16}-\frac{1}{3}\sqrt{\frac{1259}{235}},\frac{5}{16}+\frac{1}{3}\sqrt{\frac{1259}{235}}\Bigg).$$
The points $(-1,y_2),(-1,y_3)\in\Sigma_{\scriptscriptstyle L}$ and $(1,y_0),(1,y_1)\in\Sigma_{\scriptscriptstyle R}$ are crossing points because
\begin{equation*}
\begin{aligned} 
     &\langle X_{\scriptscriptstyle L}(-1,y_2),(1,0)\rangle\cdot\langle X_{\scriptscriptstyle C}(-1,y_2),(1,0) \rangle \approx 1.5518>0,\\
     &\langle X_{\scriptscriptstyle L}(-1,y_3),(1,0)\rangle\cdot\langle X_{\scriptscriptstyle C}(-1,y_3),(1,0) \rangle\approx 17.4969>0 ,\\
     &\langle X_{\scriptscriptstyle C}(1,y_0),(1,0)\rangle\cdot\langle X_{\scriptscriptstyle R}(1,y_0),(1,0) \rangle\approx 10.9315>0,\\
	 &\langle X_{\scriptscriptstyle C}(1,y_1),(1,0)\rangle\cdot\langle X_{\scriptscriptstyle R}(1,y_1),(1,0) \rangle\approx 6.9452>0.
\end{aligned}
\end{equation*}

\noindent The orbit $(x_{\scriptscriptstyle R}(t),y_{\scriptscriptstyle R}(t))$ of $X_{\scriptscriptstyle R}$, such that $(x_{\scriptscriptstyle R}(0),y_{\scriptscriptstyle R}(0))=(1,y_0)$, is given by
\begin{equation*}
\begin{aligned}
      x_{\scriptscriptstyle R}(t)&=7 \cos(2t)+\frac{31}{48}\sqrt{\frac{1259}{235}}\sin(2t)-6, \\
      y_{\scriptscriptstyle R}(t)&=\Bigg(\frac{31}{48}\sqrt{\frac{1259}{235}}-14\Bigg)\cos(2t)-\Bigg(7+\frac{31}{24}\sqrt{\frac{1259}{235}}\Bigg)\sin(2t)+14. \\
\end{aligned}
\end{equation*}
The orbit $(x_{\scriptscriptstyle C_1}(t),y_{\scriptscriptstyle C_1}(t))$ of $X_{\scriptscriptstyle C}$,  such that $(x_{\scriptscriptstyle C_1}(0),y_{\scriptscriptstyle C_1}(0))=(1,y_1)$, is given by
\begin{equation*}
\begin{aligned}
      x_{\scriptscriptstyle C_1}(t)&=\frac{2}{3}\cos(2t)+\Bigg(\frac{1}{3} - \frac{31}{48}\sqrt{\frac{1259}{235}}\Bigg)\sin(2t)+\frac{1}{3}, \\
      y_{\scriptscriptstyle C_1}(t)&=\Bigg(\frac{1}{3} - \frac{31}{48}\sqrt{\frac{1259}{235}}\Bigg)\cos(2t) - \frac{1}{3} (1 + 4 \cos(t) \sin(t)). \\
\end{aligned}
\end{equation*}
The orbit $(x_{\scriptscriptstyle L}(t),y_{\scriptscriptstyle L}(t))$ of $X_{\scriptscriptstyle L}$,  such that $(x_{\scriptscriptstyle L}(0),y_{\scriptscriptstyle L}(0))=(-1,y_2)$, is given by
\begin{equation*}
\begin{aligned}
      x_{\scriptscriptstyle L}(t)&=-8\cos(2t)-\frac{4}{3}\sqrt{\frac{1259}{235}}\sin(2t)+7, \\
      y_{\scriptscriptstyle L}(t)&=\Bigg(4-\frac{1}{3}\sqrt{\frac{1259}{235}}\Bigg)\cos(2t)+\Bigg(4+\frac{4}{3}\sqrt{\frac{1259}{235}}\Bigg)\cos(t)\sin(t)-\frac{59}{16}. \\
\end{aligned}
\end{equation*}
The orbit $(x_{\scriptscriptstyle C_2}(t),y_{\scriptscriptstyle C_2}(t))$ of $X_{\scriptscriptstyle C}$, such that $(x_{\scriptscriptstyle C_2}(0),y_{\scriptscriptstyle C_2}(0))=(-1,y_3)$, is given by
\begin{equation*}
\begin{aligned}
      x_{\scriptscriptstyle C_2}(t)&=-\frac{4}{3}\cos(2t) + \Bigg(\frac{31}{48} + \frac{1}{3}\sqrt{\frac{1259}{235}}\Bigg)\sin(2t) + \frac{1}{3}, \\
      y_{\scriptscriptstyle C_2}(t)&=\Bigg(\frac{31}{48} + \frac{1}{3}\sqrt{\frac{1259}{235}}\Bigg)\cos(2t) + \frac{1}{3}(8 \cos(t) \sin(t) - 1). \\
\end{aligned}
\end{equation*}
The fly time of the orbit $(x_{\scriptscriptstyle R}(t),y_{\scriptscriptstyle R}(t))$, from $(1,y_0)\in \Sigma_{\scriptscriptstyle R}$ to $(1,y_1)\in \Sigma_{\scriptscriptstyle R}$, is $$t_{\scriptscriptstyle R}=\dfrac{1}{2}\arctan\Bigg(\dfrac{20832 \sqrt{295865}}{25320661}\Bigg).$$ 
The fly time of the orbit $(x_{\scriptscriptstyle C_1}(t),y_{\scriptscriptstyle C_1}(t))$, from $(1,y_1)\in \Sigma_{\scriptscriptstyle R}$ to $(-1,y_2)\in \Sigma_{\scriptscriptstyle L}$, is  
$$t_{\scriptscriptstyle C_1}=\dfrac{\pi}{2}-\dfrac{1}{2}\arctan\Bigg(\dfrac{96 (10810 + \sqrt{295865})}{496001}\Bigg).$$ 
The fly time of the orbit $(x_{\scriptscriptstyle L}(t),y_{\scriptscriptstyle L}(t))$, from $(-1,y_2)\in\Sigma_{\scriptscriptstyle L}$ to $(-1,y_3)\in\Sigma_{\scriptscriptstyle L}$, is $$t_{\scriptscriptstyle L}=\dfrac{1}{2}\arctan\Bigg(\dfrac{12 \sqrt{295865}}{7201}\Bigg).$$ 
Finally, the fly time of the orbit $(x_{\scriptscriptstyle C_2}(t),y_{\scriptscriptstyle C_2}(t))$, from $(-1,y_3)\in\Sigma_{\scriptscriptstyle L}$ to $(1,y_0)\in\Sigma_{\scriptscriptstyle R}$, is 
$$t_{\scriptscriptstyle C_2}=-\dfrac{1}{2}\arctan\Bigg(\dfrac{96 (\sqrt{295865}-10810)}{496001}\Bigg).$$

Using the mathematica software (see \cite{Mathematica}), we can draw the orbits $(x_i(t),y_i(t))$ for the time $t \in [0,t_i]$, $i=R,L,C_1,C_2$, i.e. we obtain the limit cycle given in Figure \ref{fig1}~(a). The Figure \ref{fig1}~(b) has been made with the help of P5 software (see \cite{P5}), and provide the phase portrait of vector field \eqref{sys2} in this case (the symbol $\circ$ indicates an invisible singular point).

\begin{figure}[ht]
	\begin{center}		
		\begin{overpic}[width=0.95\textwidth]{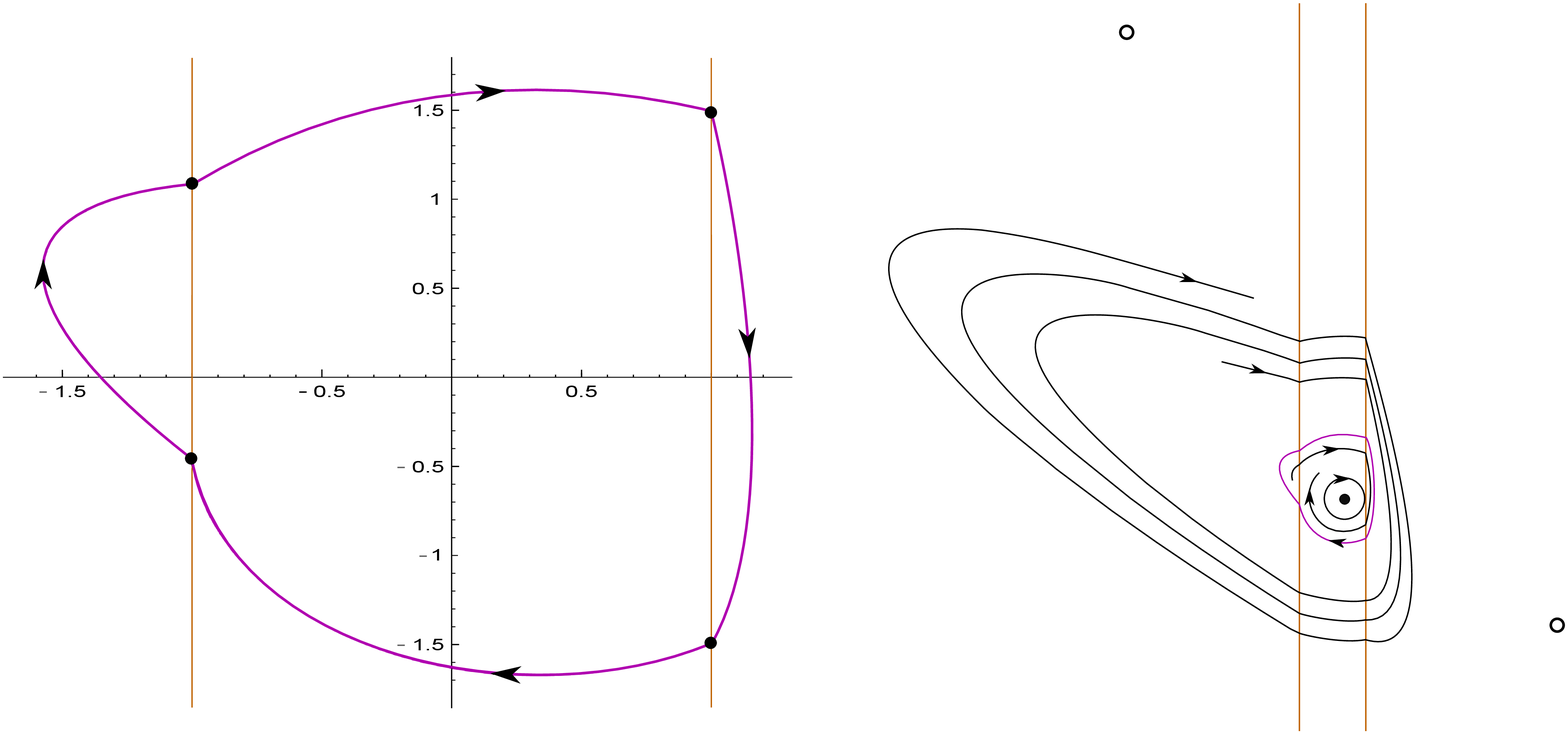}
			  	%\begin{overpic}[grid,width=0.95\textwidth]{figure_ccc4.eps}
			\put(45,44) {$\Sigma_R$}
			\put(12,44) {$\Sigma_L$}
			\put(46,38) {$(1,y_0)$}
			\put(46,4) {$(1,y_1)$}
			\put(2,15) {$(1,y_2)$}
			\put(2,36) {$(1,y_3)$}
			\put(87,47) {$\Sigma_R$}
			\put(78.5,47) {$\Sigma_L$}
			\put(26,-6) {$(a)$}
			\put(82.5,-6) {$(b)$}
		\end{overpic}
	\end{center}
	\vspace{0.7cm}
	\caption{The limit cycle of vector field (\ref{sys2}) with $a_{\scriptscriptstyle L}=4$, $b_{\scriptscriptstyle L}=8$, $\alpha_{\scriptscriptstyle L}=3/2$, $c_{\scriptscriptstyle L}=-5/2$, $\beta_{\scriptscriptstyle L}=11/4$, $a_{\scriptscriptstyle C}=0$, $b_{\scriptscriptstyle C}=2$, $\alpha_{\scriptscriptstyle C}=\beta_{\scriptscriptstyle C}=2/3$, $c_{\scriptscriptstyle C}=-2$, $a_{\scriptscriptstyle R}=4$, $b_{\scriptscriptstyle R}=2$, $c_{\scriptscriptstyle R}=-10$ and $\alpha_{\scriptscriptstyle R}=\beta_{\scriptscriptstyle R}=-4$.}
	\label{fig1}
\end{figure}  

%\begin{figure}[!h]
%\centering
%\includegraphics[scale=0.26]{figure_ccc4.pdf}
%\caption {The limit cycle of vector field (\ref{sys2}) with $a_{\scriptscriptstyle L}=4$, $b_{\scriptscriptstyle L}=8$, $\alpha_{\scriptscriptstyle L}=3/2$, $c_{\scriptscriptstyle L}=-5/2$, $\beta_{\scriptscriptstyle L}=11/4$, $a_{\scriptscriptstyle C}=0$, $b_{\scriptscriptstyle C}=2$, $\alpha_{\scriptscriptstyle C}=\beta_{\scriptscriptstyle C}=2/3$, $c_{\scriptscriptstyle C}=-2$, $a_{\scriptscriptstyle R}=4$, $b_{\scriptscriptstyle R}=2$, $c_{\scriptscriptstyle R}=-10$ and $\alpha_{\scriptscriptstyle R}=\beta_{\scriptscriptstyle R}=-4$.}\label{fig1}
%\end{figure}

\end{example}

\bigskip

\begin{example} (Case SCC)
Consider the discontinuous planar piecewise linear Hamiltonian vector field \eqref{sys2} with $a_{\scriptscriptstyle L}=b_{\scriptscriptstyle L}=1$, $\alpha_{\scriptscriptstyle L}=2/3$, $c_{\scriptscriptstyle L}=35$, $\beta_{\scriptscriptstyle L}=214/3$, $a_{\scriptscriptstyle C}=0$, $b_{\scriptscriptstyle C}=2$, $\alpha_{\scriptscriptstyle C}=\beta_{\scriptscriptstyle C}=2/3$, $c_{\scriptscriptstyle C}=-2$, $a_{\scriptscriptstyle R}=4$, $b_{\scriptscriptstyle R}=2$, $\alpha_{\scriptscriptstyle R}=\beta_{\scriptscriptstyle R}=-4$ and $c_{\scriptscriptstyle R}=-10$.
%\begin{equation}\label{ex2}
   %\begin{aligned}
      %X_{\scriptscriptstyle L}(x,y)&=\Big(x+y+\frac{2}{3},35x - y + \frac{214}{3}\Big),\quad x \leq -1, \\
      %X_{\scriptscriptstyle C}(x,y)&=\Big(2 y + \frac{2}{3},- 2 x + \frac{2}{3}\Big),\quad -1\leq x\leq 1, \\
	  %X_{\scriptscriptstyle R}(x,y)&=( 4 x + 2 y - 4, - 10 x - 4 y - 4),\quad x \geq 1.
  %\end{aligned}
%\end{equation}
The eigenvalues of the linear part of $X_i$, $i=L, C, R$,  from \eqref{sys2} for this case, are $\pm 6$, $\pm 2i$ and $\pm 2i$, respectively, i.e. we have one saddle and two centers.
%\begin{equation*}  
%\begin{aligned}
     % H_{\scriptscriptstyle L}(x,y)&=\frac{y^2}{2} -\frac{35x^2}{2} + xy + \frac{2y}{3} - \frac{214x}{3},\quad x\leq -1, \\
      %H_{\scriptscriptstyle C}(x,y)&=x^2+y^2-\frac{2x}{3} + \frac{2y}{3},\quad -1\leq x\leq 1, \\
	  %H_{\scriptscriptstyle R}(x,y)&=5x^2 + y^2 + 4xy + 4x - 4y,\quad x\geq 1,
 %\end{aligned}
%\end{equation*}

In this case, as in Example \ref{CCC}, the unique solution $(y_0,y_1,y_2,y_3)$  of system \eqref{eq5}   satisfying the condition $y_1<y_0$ and $y_2<y_3$, which is given by
$$\Bigg(\frac{2\sqrt{5}}{3},-\frac{2\sqrt{5}}{3},\frac{1-\sqrt{5}}{3},\frac{1+\sqrt{5}}{3}\Bigg),$$
correspond to the unique limit cycle of vector field \eqref{sys2}.

Note that the points $(-1,y_2),(-1,y_3)\in\Sigma_{\scriptscriptstyle L}$ and $(1,y_0),(1,y_1)\in\Sigma_{\scriptscriptstyle R}$ are crossing points. 

Now, we can computate:  the orbit $(x_{\scriptscriptstyle R}(t),y_{\scriptscriptstyle R}(t))$ of $X_{\scriptscriptstyle R}$ with $(x_{\scriptscriptstyle R}(0)$, $y_{\scriptscriptstyle R}(0))=(1,y_0)$; the orbit $(x_{\scriptscriptstyle C_1}(t),y_{\scriptscriptstyle C_1}(t))$ of $X_{\scriptscriptstyle C}$ with $(x_{\scriptscriptstyle C_1}(0),y_{\scriptscriptstyle C_1}(0))=(1,y_1)$; the orbit $(x_{\scriptscriptstyle L}(t),y_{\scriptscriptstyle L}(t))$ of $X_{\scriptscriptstyle L}$ with $(x_{\scriptscriptstyle L}(0),y_{\scriptscriptstyle L}(0))=(-1,y_2)$; and the orbit $(x_{\scriptscriptstyle C_2}(t),y_{\scriptscriptstyle C_2}(t))$ of $X_{\scriptscriptstyle C}$, with $(x_{\scriptscriptstyle C_2}(0),y_{\scriptscriptstyle C_2}(0))=(-1,y_3)$. We can also computate 
the fly times of the orbits: $(x_{\scriptscriptstyle R}(t),y_{\scriptscriptstyle R}(t))$ from $(1,y_0)\in\Sigma_{\scriptscriptstyle R}$ to $(1,y_1)\in\Sigma_{\scriptscriptstyle R}$; $(x_{\scriptscriptstyle C_1}(t),y_{\scriptscriptstyle C_1}(t))$ from $(1,y_1)\in\Sigma_{\scriptscriptstyle R}$ to $(-1,y_2)\in\Sigma_{\scriptscriptstyle L}$; $(x_{\scriptscriptstyle L}(t),y_{\scriptscriptstyle L}(t))$ from $(-1,y_2)\in\Sigma_{\scriptscriptstyle L}$ to $(-1,y_3)\in\Sigma_{\scriptscriptstyle L}$; and $(x_{\scriptscriptstyle C_2}(t),y_{\scriptscriptstyle C_2}(t))$ from $(-1,y_3)\in\Sigma_{\scriptscriptstyle L}$ to $(1,y_0)\in\Sigma_{\scriptscriptstyle R}$. Hence, using the mathematica software, we can draw the orbits $(x_i(t),y_i(t))$ for the time $t \in [0,t_i]$, $i=R,L,C_1,C_2$, i.e. we obtain the limit cycle given in Figure \ref{fig2}~(a). The Figure \ref{fig2}~(b) has been made with the help of P5 software, and provide the phase portrait of vector field \eqref{sys2} in this case.

\begin{figure}[ht]
	\begin{center}		
		\begin{overpic}[width=0.94\textwidth]{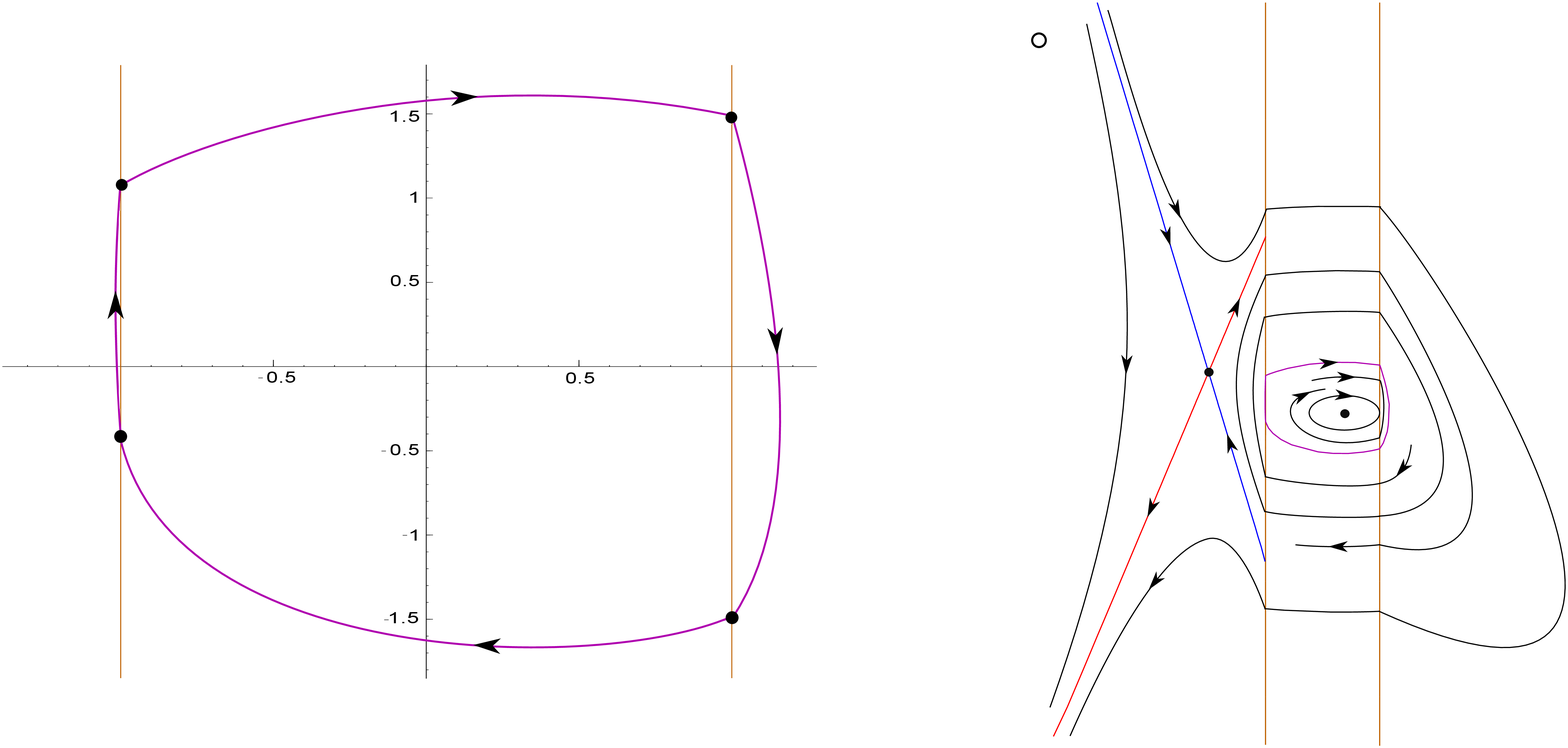}
			%\begin{overpic}[grid,width=0.95\textwidth]{figure_ccc4.eps}
			\put(46.5,45) {$\Sigma_R$}
			\put(7.5,45) {$\Sigma_L$}
			\put(48,39) {$(1,y_0)$}
			\put(48,7) {$(1,y_1)$}
			\put(-2,18) {$(1,y_2)$}
			\put(-2,36) {$(1,y_3)$}
			\put(87,49) {$\Sigma_R$}
			\put(78,49) {$\Sigma_L$}
			\put(25,-6) {$(a)$}
			\put(82.5,-6) {$(b)$}
		\end{overpic}
	\end{center}
	\vspace{0.7cm}
	\caption{The limit cycle of vector field (\ref{sys2}) with $a_{\scriptscriptstyle L}=b_{\scriptscriptstyle L}=1$, $\alpha_{\scriptscriptstyle L}=2/3$, $c_{\scriptscriptstyle L}=35$, $\beta_{\scriptscriptstyle L}=214/3$, $a_{\scriptscriptstyle C}=0$, $b_{\scriptscriptstyle C}=2$, $\alpha_{\scriptscriptstyle C}=\beta_{\scriptscriptstyle C}=2/3$, $c_{\scriptscriptstyle C}=-2$, $a_{\scriptscriptstyle R}=4$, $b_{\scriptscriptstyle R}=2$, $\alpha_{\scriptscriptstyle R}=\beta_{\scriptscriptstyle R}=-4$ and $c_{\scriptscriptstyle R}=-10$.}\label{fig2}
\end{figure}  

%\begin{figure}[!h]
%\centering
%\includegraphics[scale=0.2]{figure_scc4.pdf}
%\caption {The limit cycle of vector field (\ref{sys2}) with $a_{\scriptscriptstyle L}=b_{\scriptscriptstyle L}=1$, $\alpha_{\scriptscriptstyle L}=2/3$, $c_{\scriptscriptstyle L}=35$, $\beta_{\scriptscriptstyle L}=214/3$, $a_{\scriptscriptstyle C}=0$, $b_{\scriptscriptstyle C}=2$, $\alpha_{\scriptscriptstyle C}=\beta_{\scriptscriptstyle C}=2/3$, $c_{\scriptscriptstyle C}=-2$, $a_{\scriptscriptstyle R}=4$, $b_{\scriptscriptstyle R}=2$, $\alpha_{\scriptscriptstyle R}=\beta_{\scriptscriptstyle R}=-4$ and $c_{\scriptscriptstyle R}=-10$.}\label{fig2}
%\end{figure}

\end{example}

\bigskip

\begin{example} (Case SCS)
Consider the discontinuous planar piecewise linear Hamiltonian vector field \eqref{sys2} with $a_{\scriptscriptstyle L}=b_{\scriptscriptstyle L}=1$, $\alpha_{\scriptscriptstyle L}=3/5$, $c_{\scriptscriptstyle L}=35$, $\beta_{\scriptscriptstyle L}=357/5$, $a_{\scriptscriptstyle C}=0$, $b_{\scriptscriptstyle C}=2$, $\alpha_{\scriptscriptstyle C}=\beta_{\scriptscriptstyle C}=1$, $c_{\scriptscriptstyle C}=-2$, $a_{\scriptscriptstyle R}=b_{\scriptscriptstyle R}=1$, $\alpha_{\scriptscriptstyle R}=-1$, $c_{\scriptscriptstyle R}=15$ and $\beta_{\scriptscriptstyle R}=-31$.
%\begin{equation}\label{ex3}
 % \begin{aligned}
      %X_{\scriptscriptstyle L}(x,y)&=\Big(x+y+\frac{3}{5},35x-y+\frac{357}{5}\Big),\quad x \leq -1, \\
      %X_{\scriptscriptstyle C}(x,y)&=(2 y+1,- 2 x+1),\quad -1\leq x\leq 1, \\
	  %X_{\scriptscriptstyle R}(x,y)&=(x + y-1, 15 x - y-31),\quad x \geq 1.
  %\end{aligned} 
%\end{equation} 
The eigenvalues of the linear part of $X_i$, $i=L, C, R$,  from \eqref{sys2} for this case, are $\pm 6$, $\pm 2i$ and $\pm 4$, respectively, i.e. we have two saddles and one center.
%\begin{equation*}  
%\begin{aligned}
      %H_{\scriptscriptstyle L}(x,y)&=\frac{y^2}{2}-\frac{35x^2}{2}+xy-\frac{357 x}{5}+\frac{3 y}{5},\quad x\leq -1, \\
      %H_{\scriptscriptstyle C}(x,y)&=x^2+y^2-x+y,\quad -1\leq x\leq 1, \\
	  %H_{\scriptscriptstyle R}(x,y)&=\frac{y^2}{2}-\frac{15x^2}{2}+xy+31x-y,\quad x\geq 1,
 %\end{aligned}
%\end{equation*}

In this case, as in Example \ref{CCC}, the unique solution $(y_0,y_1,y_2,y_3)$  of system \eqref{eq5}   satisfying the condition $y_1<y_0$ and $y_2<y_3$, which is given by
$$\Bigg(\frac{18}{5}\sqrt{\frac{2}{7}},-\frac{18}{5}\sqrt{\frac{2}{7}},\frac{2}{5}-2\sqrt{\frac{2}{7}},\frac{2}{5}+2\sqrt{\frac{2}{7}}\Bigg),$$
correspond to the unique limit cycle of vector field \eqref{sys2}.

Note that the points $(-1,y_2),(-1,y_3)\in\Sigma_{\scriptscriptstyle L}$ and $(1,y_0),(1,y_1)\in\Sigma_{\scriptscriptstyle R}$ are crossing points. 

Now, we can computate:  the orbit $(x_{\scriptscriptstyle R}(t),y_{\scriptscriptstyle R}(t))$ of $X_{\scriptscriptstyle R}$ with $(x_{\scriptscriptstyle R}(0)$, $y_{\scriptscriptstyle R}(0))=(1,y_0)$; the orbit $(x_{\scriptscriptstyle C_1}(t),y_{\scriptscriptstyle C_1}(t))$ of $X_{\scriptscriptstyle C}$ with $(x_{\scriptscriptstyle C_1}(0),y_{\scriptscriptstyle C_1}(0))=(1,y_1)$; the orbit $(x_{\scriptscriptstyle L}(t),y_{\scriptscriptstyle L}(t))$ of $X_{\scriptscriptstyle L}$ with $(x_{\scriptscriptstyle L}(0),y_{\scriptscriptstyle L}(0))=(-1,y_2)$; and the orbit $(x_{\scriptscriptstyle C_2}(t),y_{\scriptscriptstyle C_2}(t))$ of $X_{\scriptscriptstyle C}$, with $(x_{\scriptscriptstyle C_2}(0),y_{\scriptscriptstyle C_2}(0))=(-1,y_3)$. We can also computate 
the fly times of the orbits: $(x_{\scriptscriptstyle R}(t),y_{\scriptscriptstyle R}(t))$ from $(1,y_0)\in\Sigma_{\scriptscriptstyle R}$ to $(1,y_1)\in\Sigma_{\scriptscriptstyle R}$; $(x_{\scriptscriptstyle C_1}(t),y_{\scriptscriptstyle C_1}(t))$ from $(1,y_1)\in\Sigma_{\scriptscriptstyle R}$ to $(-1,y_2)\in\Sigma_{\scriptscriptstyle L}$; $(x_{\scriptscriptstyle L}(t),y_{\scriptscriptstyle L}(t))$ from $(-1,y_2)\in\Sigma_{\scriptscriptstyle L}$ to $(-1,y_3)\in\Sigma_{\scriptscriptstyle L}$; and $(x_{\scriptscriptstyle C_2}(t),y_{\scriptscriptstyle C_2}(t))$ from $(-1,y_3)\in\Sigma_{\scriptscriptstyle L}$ to $(1,y_0)\in\Sigma_{\scriptscriptstyle R}$. Hence, using the mathematica software, we can draw the orbits $(x_i(t),y_i(t))$ for the time $t \in [0,t_i]$, $i=R,L,C_1,C_2$, i.e. we obtain the limit cycle given in Figure \ref{fig3}~(a). The Figure \ref{fig3}~(b) has been made with the help of P5 software, and provide the phase portrait of vector field \eqref{sys2} in this case.

\begin{figure}[ht]
	\begin{center}		
		\begin{overpic}[width=0.94\textwidth]{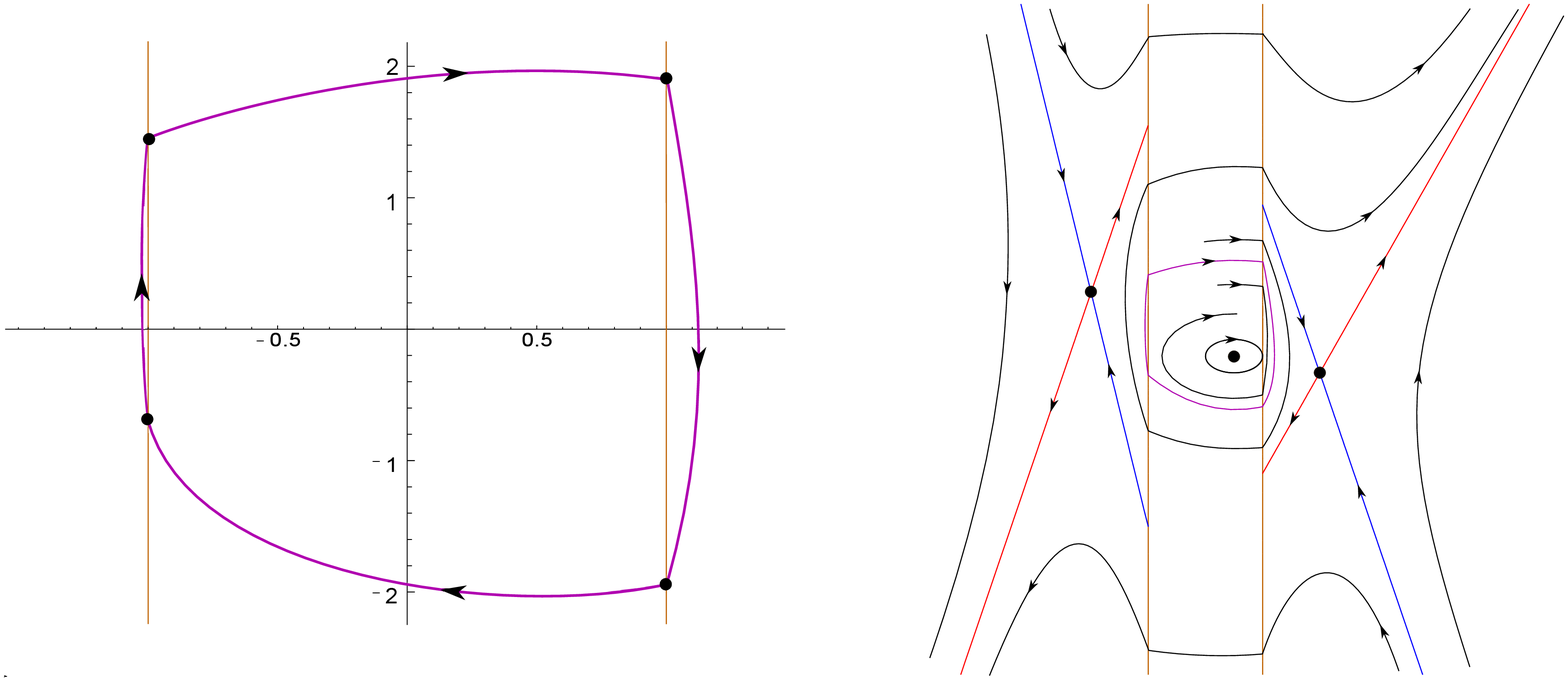}
			%\begin{overpic}[grid,width=0.95\textwidth]{figure_ccc4.eps}
			\put(42,42) {$\Sigma_R$}
			\put(9,42) {$\Sigma_L$}
			\put(44,37) {$(1,y_0)$}
			\put(44,6) {$(1,y_1)$}
			\put(-1,15) {$(1,y_2)$}
			\put(-1,34) {$(1,y_3)$}
			\put(79,45) {$\Sigma_R$}
			\put(71,45) {$\Sigma_L$}
			\put(24,-6) {$(a)$}
			\put(75,-6) {$(b)$}
		\end{overpic}
	\end{center}
	\vspace{0.7cm}
	\caption {The limit cycle of vector field (\ref{sys2}) with $a_{\scriptscriptstyle L}=b_{\scriptscriptstyle L}=1$, $\alpha_{\scriptscriptstyle L}=3/5$, $c_{\scriptscriptstyle L}=35$, $\beta_{\scriptscriptstyle L}=357/5$, $a_{\scriptscriptstyle C}=0$, $b_{\scriptscriptstyle C}=2$, $\alpha_{\scriptscriptstyle C}=\beta_{\scriptscriptstyle C}=1$, $c_{\scriptscriptstyle C}=-2$, $a_{\scriptscriptstyle R}=b_{\scriptscriptstyle R}=1$, $\alpha_{\scriptscriptstyle R}=-1$, $c_{\scriptscriptstyle R}=15$ and $\beta_{\scriptscriptstyle R}=-31$.}\label{fig3}
\end{figure}  

%\begin{figure}[!h]
%\centering
%\includegraphics[scale=0.35]{figure_scs4.pdf}
%\caption {The limit cycle of vector field (\ref{sys2}) with $a_{\scriptscriptstyle L}=b_{\scriptscriptstyle L}=1$, $\alpha_{\scriptscriptstyle L}=3/5$, $c_{\scriptscriptstyle L}=35$, $\beta_{\scriptscriptstyle L}=357/5$, $a_{\scriptscriptstyle C}=0$, $b_{\scriptscriptstyle C}=2$, $\alpha_{\scriptscriptstyle C}=\beta_{\scriptscriptstyle C}=1$, $c_{\scriptscriptstyle C}=-2$, $a_{\scriptscriptstyle R}=b_{\scriptscriptstyle R}=1$, $\alpha_{\scriptscriptstyle R}=-1$, $c_{\scriptscriptstyle R}=15$ and $\beta_{\scriptscriptstyle R}=-31$.}\label{fig3}
%\end{figure}

\end{example}

\begin{example} (Case CSC)
Consider the discontinuous planar piecewise linear Hamiltonian vector field \eqref{sys2} with $a_{\scriptscriptstyle L}=4$, $b_{\scriptscriptstyle L}=8$, $\alpha_{\scriptscriptstyle L}=2$, $c_{\scriptscriptstyle L}=-5/2$, $\beta_{\scriptscriptstyle L}=5/2$, $a_{\scriptscriptstyle C}=2/5$, $b_{\scriptscriptstyle C}=24/5$, $\alpha_{\scriptscriptstyle C}=-9/5$, $c_{\scriptscriptstyle C}=4/5$, $\beta_{\scriptscriptstyle C}=-4/15$, $a_{\scriptscriptstyle R}=8$, $b_{\scriptscriptstyle R}=10$ and $\alpha_{\scriptscriptstyle R}=c_{\scriptscriptstyle R}=\beta_{\scriptscriptstyle R}=-8$. 
%\begin{equation}\label{ex4}
  %\begin{aligned}
      %X_{\scriptscriptstyle L}(x,y)&=\Big(4x+8y+2,-\frac{5x}{2}-4y+\frac{5}{2}\Big),\quad x \leq -1, \\
      %X_{\scriptscriptstyle C}(x,y)&=\Big(\frac{2x}{5}+\frac{24y}{5}-\frac{9}{5},\frac{4x}{5}-\frac{2y}{5}-\frac{4}{15}\Big),\quad -1\leq x\leq 1, \\
	  %X_{\scriptscriptstyle R}(x,y)&=( 8 x + 10 y-8,- 8 x - 8 y - 8),\quad x \geq 1.
  %\end{aligned}
%\end{equation} 
The eigenvalues of the linear part of $X_i$, $i=L, C, R$,  from \eqref{sys2} for this case, are $\pm 2i$, $\pm 2$ and $\pm 4i$, respectively, i.e. we have one saddle and two centers.
%\begin{equation}  
%\begin{aligned}
      %H_{\scriptscriptstyle L}(x,y)&=\frac{5x^2}{4}+4y^2+4xy-\frac{5x}{2}+2y,\quad x\leq -1, \\
      %H_{\scriptscriptstyle C}(x,y)&=\frac{12y^2}{5}-\frac{2x^2}{5}+\frac{2xy}{5}+\frac{4x}{15}-\frac{9y}{5},\quad -1\leq x\leq 1, \\
	 % H_{\scriptscriptstyle R}(x,y)&=4 x^2 + 5 y^2 + 8 x y +8 x- 8 y,\quad x\geq 1,
 %\end{aligned}
%\end{equation}

In this case, as in Example \ref{CCC}, the unique solution $(y_0,y_1,y_2,y_3)$  of system \eqref{eq5}   satisfying the condition $y_1<y_0$ and $y_2<y_3$, which is given by
$$\Bigg(\frac{5}{12}\sqrt{\frac{7}{3}},-\frac{5}{12}\sqrt{\frac{7}{3}},\frac{1}{4}-\frac{7\sqrt{21}}{36},\frac{1}{4}+\frac{7\sqrt{21}}{36}\Bigg),$$
correspond to the unique limit cycle of vector field \eqref{sys2}.

Note that the points $(-1,y_2),(-1,y_3)\in\Sigma_{\scriptscriptstyle L}$ and $(1,y_0),(1,y_1)\in\Sigma_{\scriptscriptstyle R}$ are crossing points.

Now, we can computate:  the orbit $(x_{\scriptscriptstyle R}(t),y_{\scriptscriptstyle R}(t))$ of $X_{\scriptscriptstyle R}$ with $(x_{\scriptscriptstyle R}(0)$, $y_{\scriptscriptstyle R}(0))=(1,y_0)$; the orbit $(x_{\scriptscriptstyle C_1}(t),y_{\scriptscriptstyle C_1}(t))$ of $X_{\scriptscriptstyle C}$ with $(x_{\scriptscriptstyle C_1}(0),y_{\scriptscriptstyle C_1}(0))=(1,y_1)$; the orbit $(x_{\scriptscriptstyle L}(t),y_{\scriptscriptstyle L}(t))$ of $X_{\scriptscriptstyle L}$ with $(x_{\scriptscriptstyle L}(0),y_{\scriptscriptstyle L}(0))=(-1,y_2)$; and the orbit $(x_{\scriptscriptstyle C_2}(t),y_{\scriptscriptstyle C_2}(t))$ of $X_{\scriptscriptstyle C}$, with $(x_{\scriptscriptstyle C_2}(0),y_{\scriptscriptstyle C_2}(0))=(-1,y_3)$. We can also computate the fly times of the orbits: $(x_{\scriptscriptstyle R}(t),y_{\scriptscriptstyle R}(t))$ from $(1,y_0)\in\Sigma_{\scriptscriptstyle R}$ to $(1,y_1)\in\Sigma_{\scriptscriptstyle R}$; $(x_{\scriptscriptstyle C_1}(t),y_{\scriptscriptstyle C_1}(t))$ from $(1,y_1)\in\Sigma_{\scriptscriptstyle R}$ to $(-1,y_2)\in\Sigma_{\scriptscriptstyle L}$; $(x_{\scriptscriptstyle L}(t),y_{\scriptscriptstyle L}(t))$ from $(-1,y_2)\in\Sigma_{\scriptscriptstyle L}$ to $(-1,y_3)\in\Sigma_{\scriptscriptstyle L}$; and $(x_{\scriptscriptstyle C_2}(t),y_{\scriptscriptstyle C_2}(t))$ from $(-1,y_3)\in\Sigma_{\scriptscriptstyle L}$ to $(1,y_0)\in\Sigma_{\scriptscriptstyle R}$. Hence, using the mathematica software, we can draw the orbits $(x_i(t),y_i(t))$ for the time $t \in [0,t_i]$, $i=R,L,C_1,C_2$, i.e. we obtain the limit cycle given in Figure \ref{fig4}~(a). The Figure \ref{fig4}~(b) has been made with the help of P5 software, and provide the phase portrait of vector field \eqref{sys2} in this case.

\bigskip

\begin{figure}[ht]
	\begin{center}		
		\begin{overpic}[width=0.94\textwidth]{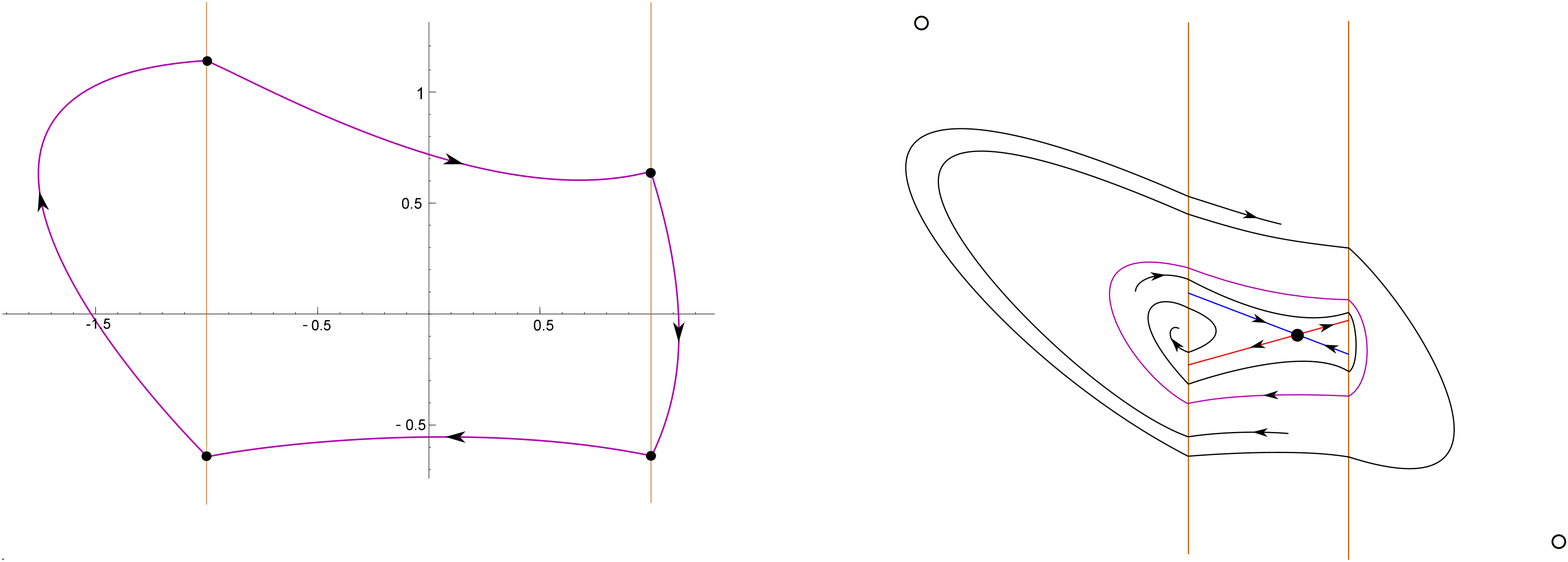}
			%\begin{overpic}[grid,width=0.95\textwidth]{figure_ccc4.eps}
			\put(41,37) {$\Sigma_R$}
			\put(13,37) {$\Sigma_L$}
			\put(43,24) {$(1,y_0)$}
			\put(43,6) {$(1,y_1)$}
			\put(3,5) {$(1,y_2)$}
			\put(3,33) {$(1,y_3)$}
			\put(84,36) {$\Sigma_R$}
			\put(73,36) {$\Sigma_L$}
			\put(25,-6) {$(a)$}
			\put(79,-6) {$(b)$}
		\end{overpic}
	\end{center}
	\vspace{0.7cm}
\caption {The limit cycle of vector field (\ref{sys2}) with $a_{\scriptscriptstyle L}=4$, $b_{\scriptscriptstyle L}=8$, $\alpha_{\scriptscriptstyle L}=2$, $c_{\scriptscriptstyle L}=-5/2$, $\beta_{\scriptscriptstyle L}=5/2$, $a_{\scriptscriptstyle C}=2/5$, $b_{\scriptscriptstyle C}=24/5$, $\alpha_{\scriptscriptstyle C}=-9/5$, $c_{\scriptscriptstyle C}=4/5$, $\beta_{\scriptscriptstyle C}=-4/15$, $a_{\scriptscriptstyle R}=8$, $b_{\scriptscriptstyle R}=10$ and $\alpha_{\scriptscriptstyle R}=c_{\scriptscriptstyle R}=\beta_{\scriptscriptstyle R}=-8$.}\label{fig4}
\end{figure}  

%\begin{figure}[!h]
%\centering
%\includegraphics[scale=0.15]{figure_csc4.pdf}
%\caption {The limit cycle of vector field (\ref{sys2}) with $a_{\scriptscriptstyle L}=4$, $b_{\scriptscriptstyle L}=8$, $\alpha_{\scriptscriptstyle L}=2$, $c_{\scriptscriptstyle L}=-5/2$, $\beta_{\scriptscriptstyle L}=5/2$, $a_{\scriptscriptstyle C}=2/5$, $b_{\scriptscriptstyle C}=24/5$, $\alpha_{\scriptscriptstyle C}=-9/5$, $c_{\scriptscriptstyle C}=4/5$, $\beta_{\scriptscriptstyle C}=-4/15$, $a_{\scriptscriptstyle R}=8$, $b_{\scriptscriptstyle R}=10$ and $\alpha_{\scriptscriptstyle R}=c_{\scriptscriptstyle R}=\beta_{\scriptscriptstyle R}=-8$.}\label{fig4}
%\end{figure}

\end{example}

\bigskip

\begin{example} (Case SSS)
Consider the discontinuous planar piecewise linear Hamiltonian vector field \eqref{sys2} with $a_{\scriptscriptstyle L}=\alpha_{\scriptscriptstyle L}=-2/3$, $b_{\scriptscriptstyle L}=4/3$, $c_{\scriptscriptstyle L}=8/3$, $\beta_{\scriptscriptstyle L}=35/3$, $a_{\scriptscriptstyle C}=2/11$, $b_{\scriptscriptstyle C}=120/11$, $\alpha_{\scriptscriptstyle C}=-41/11$, $c_{\scriptscriptstyle C}=4/11$, $\beta_{\scriptscriptstyle C}=-4/33$, $a_{\scriptscriptstyle R}=-2/11$, $b_{\scriptscriptstyle R}=4/11$, $\alpha_{\scriptscriptstyle R}=1/5$, $c_{\scriptscriptstyle R}=120/11$ and $\beta_{\scriptscriptstyle R}=-749/55$. 
%\begin{equation}\label{ex5}
 % \begin{aligned}
      %X_{\scriptscriptstyle L}(x,y)&=\Big(-\frac{2x}{3}+\frac{4y}{3}-\frac{2}{3},\frac{8x}{3}+\frac{2y}{3}+\frac{35}{3}\Big),\quad x \leq -1, \\
      %X_{\scriptscriptstyle C}(x,y)&=\Big(\frac{2x}{11}+\frac{120y}{11}-\frac{41}{11},\frac{4x}{11}-\frac{2y}{11}-\frac{4}{33}\Big),\quad -1\leq x\leq 1, \\
	  %X_{\scriptscriptstyle R}(x,y)&=\Big(-\frac{2x}{11}+\frac{4y}{11}+\frac{1}{5},\frac{120x}{11}+\frac{2y}{11}-\frac{749}{55}\Big),\quad x \geq 1.
  %\end{aligned}
%\end{equation} 
The eigenvalues of the linear part of $X_i$, $i=L, C, R$,  from \eqref{sys2} for this case, are $\pm 2$, $\pm 2$ and $\pm 2$, respectively, i.e. we have three saddles.
%\begin{equation*}  
%\begin{aligned}
      %H_{\scriptscriptstyle L}(x,y)&=\frac{2y^2}{3}-\frac{4x^2}{3}-\frac{2xy}{3}-\frac{35x}{3}-\frac{2y}{3},\quad x\leq -1, \\
      %H_{\scriptscriptstyle C}(x,y)&=\frac{60y^2}{11}-\frac{2x^2}{11}+\frac{2xy}{11}+\frac{4x}{33}-\frac{41y}{11},\quad -1\leq x\leq 1, \\
	  %H_{\scriptscriptstyle R}(x,y)&=\frac{2y^2}{11}-\frac{60x^2}{11}-\frac{2xy}{11}+\frac{749x}{55}+\frac{y}{5},\quad x\geq 1,
 %\end{aligned}
%\end{equation*}
In this case, as in Example \ref{CCC}, the unique solution $(y_0,y_1,y_2,y_3)$  of system \eqref{eq5} satisfying the condition $y_1<y_0$ and $y_2<y_3$, which is given by
$$\Bigg(\frac{43\sqrt{26}-12}{240},-\frac{43\sqrt{26}+12}{240},-\frac{3}{8} \sqrt{\frac{13}{2}},\frac{3}{8} \sqrt{\frac{13}{2}}\Bigg),$$
correspond to the unique limit cycle of vector field \eqref{sys2}.

Note that the points $(-1,y_2),(-1,y_3)\in\Sigma_{\scriptscriptstyle L}$ and $(1,y_0),(1,y_1)\in\Sigma_{\scriptscriptstyle R}$ are crossing points.

Now, we can computate:  the orbit $(x_{\scriptscriptstyle R}(t),y_{\scriptscriptstyle R}(t))$ of $X_{\scriptscriptstyle R}$ with $(x_{\scriptscriptstyle R}(0)$, $y_{\scriptscriptstyle R}(0))=(1,y_0)$; the orbit $(x_{\scriptscriptstyle C_1}(t),y_{\scriptscriptstyle C_1}(t))$ of $X_{\scriptscriptstyle C}$ with $(x_{\scriptscriptstyle C_1}(0),y_{\scriptscriptstyle C_1}(0))=(1,y_1)$; the orbit $(x_{\scriptscriptstyle L}(t),y_{\scriptscriptstyle L}(t))$ of $X_{\scriptscriptstyle L}$ with $(x_{\scriptscriptstyle L}(0),y_{\scriptscriptstyle L}(0))=(-1,y_2)$; and the orbit $(x_{\scriptscriptstyle C_2}(t),y_{\scriptscriptstyle C_2}(t))$ of $X_{\scriptscriptstyle C}$, with $(x_{\scriptscriptstyle C_2}(0),y_{\scriptscriptstyle C_2}(0))=(-1,y_3)$. We can also computate the fly times of the orbits: $(x_{\scriptscriptstyle R}(t),y_{\scriptscriptstyle R}(t))$ from $(1,y_0)\in\Sigma_{\scriptscriptstyle R}$ to $(1,y_1)\in\Sigma_{\scriptscriptstyle R}$; $(x_{\scriptscriptstyle C_1}(t),y_{\scriptscriptstyle C_1}(t))$ from $(1,y_1)\in\Sigma_{\scriptscriptstyle R}$ to $(-1,y_2)\in\Sigma_{\scriptscriptstyle L}$; $(x_{\scriptscriptstyle L}(t),y_{\scriptscriptstyle L}(t))$ from $(-1,y_2)\in\Sigma_{\scriptscriptstyle L}$ to $(-1,y_3)\in\Sigma_{\scriptscriptstyle L}$; and $(x_{\scriptscriptstyle C_2}(t),y_{\scriptscriptstyle C_2}(t))$ from $(-1,y_3)\in\Sigma_{\scriptscriptstyle L}$ to $(1,y_0)\in\Sigma_{\scriptscriptstyle R}$. Hence, using the mathematica software, we can draw the orbits $(x_i(t),y_i(t))$ for the time $t \in [0,t_i]$, $i=R,L,C_1,C_2$, i.e. we obtain the limit cycle given in Figure \ref{fig5}~(a). The Figure \ref{fig5}~(b) has been made with the help of P5 software, and provide the phase portrait of vector field \eqref{sys2} in this case..

\begin{figure}[ht]
	\begin{center}		
		\begin{overpic}[width=0.9\textwidth]{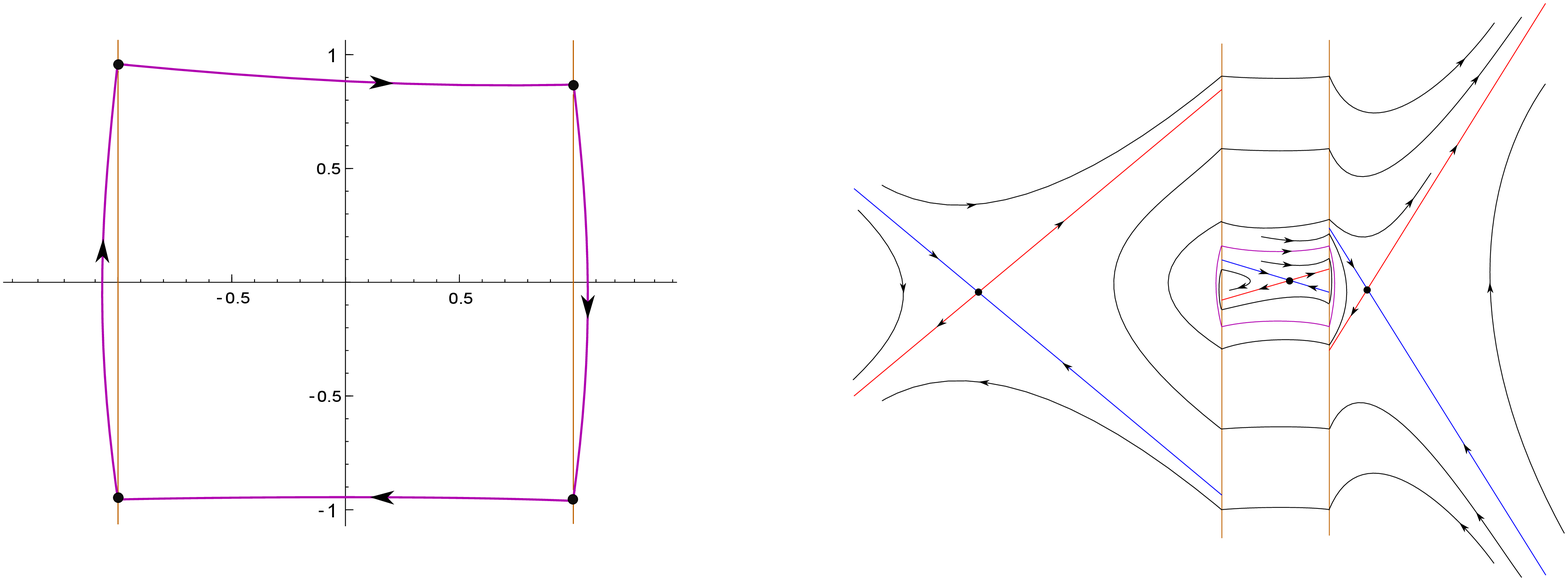}
			%\begin{overpic}[grid,width=0.95\textwidth]{figure_ccc4.eps}
			\put(36,36) {$\Sigma_R$}
			\put(7,36) {$\Sigma_L$}
			\put(38,31) {$(1,y_0)$}
			\put(38,4) {$(1,y_1)$}
			\put(-3,4) {$(1,y_2)$}
			\put(-3,33) {$(1,y_3)$}
			\put(83,36) {$\Sigma_R$}
			\put(76,36) {$\Sigma_L$}
			\put(20,-6) {$(a)$}
			\put(80,-6) {$(b)$}
		\end{overpic}
	\end{center}
	\vspace{0.7cm}
	\caption {The limit cycle of vector field (\ref{sys2}) with $a_{\scriptscriptstyle L}=\alpha_{\scriptscriptstyle L}=-2/3$, $b_{\scriptscriptstyle L}=4/3$, $c_{\scriptscriptstyle L}=8/3$, $\beta_{\scriptscriptstyle L}=35/3$, $a_{\scriptscriptstyle C}=2/11$, $b_{\scriptscriptstyle C}=120/11$, $\alpha_{\scriptscriptstyle C}=-41/11$, $c_{\scriptscriptstyle C}=4/11$, $\beta_{\scriptscriptstyle C}=-4/33$, $a_{\scriptscriptstyle R}=-2/11$, $b_{\scriptscriptstyle R}=4/11$, $\alpha_{\scriptscriptstyle R}=1/5$, $c_{\scriptscriptstyle R}=120/11$ and $\beta_{\scriptscriptstyle R}=-749/55$.}\label{fig5}
\end{figure} 

%\begin{figure}[!h]
%\centering
%\includegraphics[scale=0.25]{figure_sss4.pdf}
%\caption {The limit cycle of vector field (\ref{sys2}) with $a_{\scriptscriptstyle L}=\alpha_{\scriptscriptstyle L}=-2/3$, $b_{\scriptscriptstyle L}=4/3$, $c_{\scriptscriptstyle L}=8/3$, $\beta_{\scriptscriptstyle L}=35/3$, $a_{\scriptscriptstyle C}=2/11$, $b_{\scriptscriptstyle C}=120/11$, $\alpha_{\scriptscriptstyle C}=-41/11$, $c_{\scriptscriptstyle C}=4/11$, $\beta_{\scriptscriptstyle C}=-4/33$, $a_{\scriptscriptstyle R}=-2/11$, $b_{\scriptscriptstyle R}=4/11$, $\alpha_{\scriptscriptstyle R}=1/5$, $c_{\scriptscriptstyle R}=120/11$ and $\beta_{\scriptscriptstyle R}=-749/55$.}\label{fig5}
%\end{figure}

\end{example}

\bigskip

\begin{example} (Case SSC)
Consider the discontinuous planar piecewise linear Hamiltonian vector field \eqref{sys2} with $a_{\scriptscriptstyle L}=\alpha_{\scriptscriptstyle L}=-2/3$, $b_{\scriptscriptstyle L}=4/3$, $c_{\scriptscriptstyle L}=8/3$, $\beta_{\scriptscriptstyle L}=35/3$, $a_{\scriptscriptstyle C}=2/11$, $b_{\scriptscriptstyle C}=120/11$, $\alpha_{\scriptscriptstyle C}=-41/11$, $c_{\scriptscriptstyle C}=4/11$, $\beta_{\scriptscriptstyle C}=-4/33$, $a_{\scriptscriptstyle R}=8$, $b_{\scriptscriptstyle R}=10$, $\alpha_{\scriptscriptstyle R}=-7$ and $c_{\scriptscriptstyle R}=\beta_{\scriptscriptstyle R}=-8$.
%\begin{equation}\label{ex6}
  %\begin{aligned}
      %X_{\scriptscriptstyle L}(x,y)&=\Big(-\frac{2x}{3}+\frac{4y}{3}-\frac{2}{3},\frac{8x}{3}+\frac{2y}{3}+\frac{35}{3}\Big),\quad x \leq -1, \\
      %X_{\scriptscriptstyle C}(x,y)&=\Big(\frac{2x}{11}+\frac{120y}{11}-\frac{41}{11},\frac{4x}{11}-\frac{2y}{11}-\frac{4}{33}\Big),\quad -1\leq x\leq 1, \\
	  %X_{\scriptscriptstyle R}(x,y)&=(8x+10y-7,-8x-8y-8),\quad x \geq 1.
  %\end{aligned}
%\end{equation} 
The eigenvalues of the linear part of $X_i$, $i=L, C, R$,  from \eqref{sys2} for this case, are $\pm 2$, $\pm 2$ and $\pm 4i$, respectively, i.e. we have two saddles and one center.
%\begin{equation*}  
%\begin{aligned}
      %H_{\scriptscriptstyle L}(x,y)&=\frac{2y^2}{3}-\frac{4x^2}{3}-\frac{2xy}{3}-\frac{35x}{3}-\frac{2y}{3},\quad x\leq -1, \\
      %H_{\scriptscriptstyle C}(x,y)&=\frac{60y^2}{11}-\frac{2x^2}{11}+\frac{2xy}{11}+\frac{4x}{33}-\frac{41y}{11},\quad -1\leq x\leq 1, \\
	  %H_{\scriptscriptstyle R}(x,y)&=5y^2+4x^2+8xy+8x-7y,\quad x\geq 1,
 %\end{aligned}
%\end{equation*}
In this case, as in Example \ref{CCC}, the unique solution $(y_0,y_1,y_2,y_3)$  of system \eqref{eq5} satisfying the condition $y_1<y_0$ and $y_2<y_3$, which is given by
$$\Bigg(\frac{43}{24}\sqrt{\frac{43}{470}}-\frac{1}{10},-\frac{43}{24}\sqrt{\frac{43}{470}}-\frac{1}{10},-\frac{17}{8} \sqrt{\frac{43}{470}},\frac{17}{8} \sqrt{\frac{43}{470}}\Bigg),$$
correspond to the unique limit cycle of vector field \eqref{sys2}.

Note that the points $(-1,y_2),(-1,y_3)\in\Sigma_{\scriptscriptstyle L}$ and $(1,y_0),(1,y_1)\in\Sigma_{\scriptscriptstyle R}$ are crossing points.

Now, we can computate:  the orbit $(x_{\scriptscriptstyle R}(t),y_{\scriptscriptstyle R}(t))$ of $X_{\scriptscriptstyle R}$ with $(x_{\scriptscriptstyle R}(0)$, $y_{\scriptscriptstyle R}(0))=(1,y_0)$; the orbit $(x_{\scriptscriptstyle C_1}(t),y_{\scriptscriptstyle C_1}(t))$ of $X_{\scriptscriptstyle C}$ with $(x_{\scriptscriptstyle C_1}(0),y_{\scriptscriptstyle C_1}(0))=(1,y_1)$; the orbit $(x_{\scriptscriptstyle L}(t),y_{\scriptscriptstyle L}(t))$ of $X_{\scriptscriptstyle L}$ with $(x_{\scriptscriptstyle L}(0),y_{\scriptscriptstyle L}(0))=(-1,y_2)$; and the orbit $(x_{\scriptscriptstyle C_2}(t),y_{\scriptscriptstyle C_2}(t))$ of $X_{\scriptscriptstyle C}$, with $(x_{\scriptscriptstyle C_2}(0),y_{\scriptscriptstyle C_2}(0))=(-1,y_3)$. We can also computate the fly times of the orbits: $(x_{\scriptscriptstyle R}(t),y_{\scriptscriptstyle R}(t))$ from $(1,y_0)\in\Sigma_{\scriptscriptstyle R}$ to $(1,y_1)\in\Sigma_{\scriptscriptstyle R}$; $(x_{\scriptscriptstyle C_1}(t),y_{\scriptscriptstyle C_1}(t))$ from $(1,y_1)\in\Sigma_{\scriptscriptstyle R}$ to $(-1,y_2)\in\Sigma_{\scriptscriptstyle L}$; $(x_{\scriptscriptstyle L}(t),y_{\scriptscriptstyle L}(t))$ from $(-1,y_2)\in\Sigma_{\scriptscriptstyle L}$ to $(-1,y_3)\in\Sigma_{\scriptscriptstyle L}$; and $(x_{\scriptscriptstyle C_2}(t),y_{\scriptscriptstyle C_2}(t))$ from $(-1,y_3)\in\Sigma_{\scriptscriptstyle L}$ to $(1,y_0)\in\Sigma_{\scriptscriptstyle R}$. Hence, using the mathematica software, we can draw the orbits $(x_i(t),y_i(t))$ for the time $t \in [0,t_i]$, $i=R,L,C_1,C_2$, i.e. we obtain the limit cycle given in Figure \ref{fig6}~(a). The Figure \ref{fig6}~(b) has been made with the help of P5 software, and provide the phase portrait of vector field \eqref{sys2} in this case.

\begin{figure}[ht]
	\begin{center}		
		\begin{overpic}[width=0.9\textwidth]{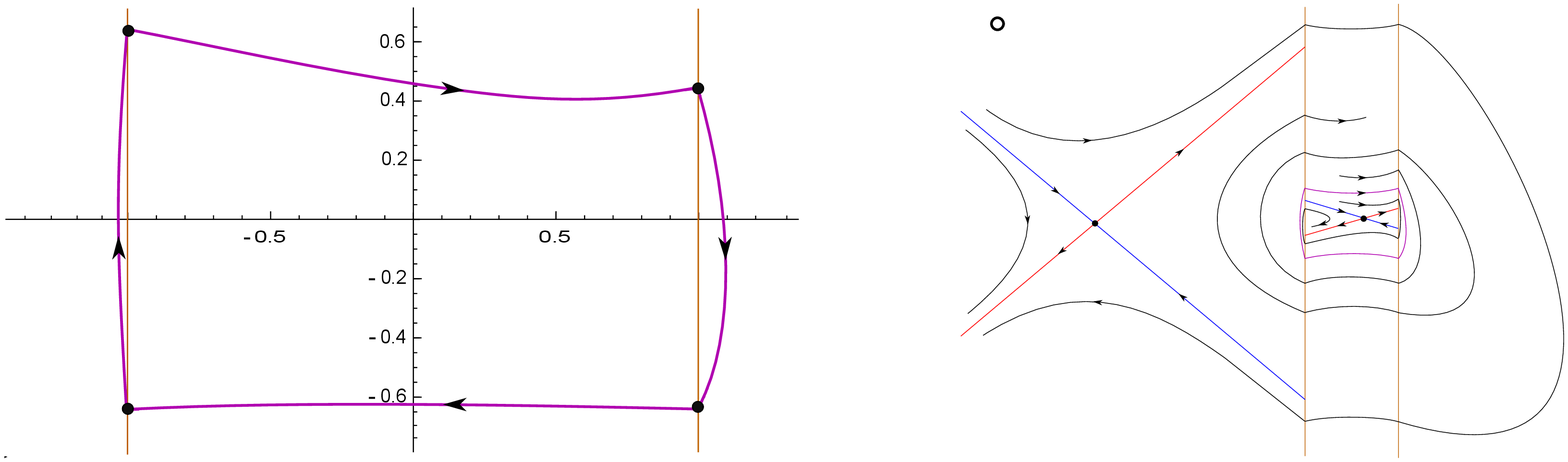}
			%\begin{overpic}[grid,width=0.95\textwidth]{figure_ccc4.eps}
			\put(43,30) {$\Sigma_R$}
			\put(7,30) {$\Sigma_L$}
			\put(46,23) {$(1,y_0)$}
			\put(46,3) {$(1,y_1)$}
			\put(-2,3) {$(1,y_2)$}
			\put(-2,26) {$(1,y_3)$}
			\put(89,31) {$\Sigma_R$}
			\put(79,31) {$\Sigma_L$}
			\put(24,-6) {$(a)$}
			\put(84,-6) {$(b)$}
		\end{overpic}
	\end{center}
	\vspace{0.7cm}
\caption {The limit cycle of vector field (\ref{sys2}) with $a_{\scriptscriptstyle L}=\alpha_{\scriptscriptstyle L}=-2/3$, $b_{\scriptscriptstyle L}=4/3$, $c_{\scriptscriptstyle L}=8/3$, $\beta_{\scriptscriptstyle L}=35/3$, $a_{\scriptscriptstyle C}=2/11$, $b_{\scriptscriptstyle C}=120/11$, $\alpha_{\scriptscriptstyle C}=-41/11$, $c_{\scriptscriptstyle C}=4/11$, $\beta_{\scriptscriptstyle C}=-4/33$, $a_{\scriptscriptstyle R}=8$, $b_{\scriptscriptstyle R}=10$, $\alpha_{\scriptscriptstyle R}=-7$ and $c_{\scriptscriptstyle R}=\beta_{\scriptscriptstyle R}=-8$.}\label{fig6}
\end{figure} 

%\begin{figure}[!h]
%\centering
%\includegraphics[scale=0.35]{figure_ssc4.pdf}
%\caption {The limit cycle of vector field (\ref{sys2}) with $a_{\scriptscriptstyle L}=\alpha_{\scriptscriptstyle L}=-2/3$, $b_{\scriptscriptstyle L}=4/3$, $c_{\scriptscriptstyle L}=8/3$, $\beta_{\scriptscriptstyle L}=35/3$, $a_{\scriptscriptstyle C}=2/11$, $b_{\scriptscriptstyle C}=120/11$, $\alpha_{\scriptscriptstyle C}=-41/11$, $c_{\scriptscriptstyle C}=4/11$, $\beta_{\scriptscriptstyle C}=-4/33$, $a_{\scriptscriptstyle R}=8$, $b_{\scriptscriptstyle R}=10$, $\alpha_{\scriptscriptstyle R}=-7$ and $c_{\scriptscriptstyle R}=\beta_{\scriptscriptstyle R}=-8$.}\label{fig6}
%\end{figure}

\end{example}

\medskip

\section{Acknowledgments}

The first author is partially supported by S\~ao Paulo Research Foundation (FAPESP) grants 19/10269-3 and 18/19726-5.
The second author is supported by CAPES grants 88882.434343/2019-01.

%\addcontentsline{toc}{chapter}{Bibliografia}
%\bibliographystyle{siam}
%\bibliography{cent_discontinuous.bib}

\addcontentsline{toc}{chapter}{Bibliografia}
\begin{thebibliography}{80}

\bibitem{And} Andronov, A., Vitt, A. and Khaikin, S.: {\it Theory of Oscillations}, Pergamon Press, Oxford (1966).

\bibitem{Bra} Braga, D. and Mello, L.: {\it Limit cycles in a family of discontinuous piecewise linear differential systems with two zones in the plane}, Nonlin. Dyn. 73, 1283-1288 (2013).

\bibitem{Bra2} Braga, D. and Mello, L.: {\it More than three limit cycles in discontinuous piecewise linear differential systems with two zones in the plane}, Int. J. Bifur. Chaos Appl. Sci. Eng. 24, 1450056, 10 pp (2014).

\bibitem{Buz} Buzzi, C., Pessoa, C. and Torregrosa, J.: {\it Piecewise linear perturbations of a linear center}, Discrete Contin. Dyn. Syst. 33, 3915-3936 (2013).

\bibitem{Cas} Castillo, J., Llibre, J. and Verduzco, F.: {\it The pseudo-Hopf bifurcation for planar discontinuous piecewise 
	ar differential systems}, Nonlin. Dyn. 90, 1829-1840 (2017).

\bibitem{di} di Bernardo, M., Budd, C. J., Champneys, A. R. and Kowalczyk, P.: {\it Piecewise-Smooth Dynamical Systems: Theory and Applications}, Springer (2008).

\bibitem{Mel} Fonseca, A., Llibre, J. and Mello, L.: {\it Limit cycles in planar piecewise linear Hamiltonian systems with three zones without equilibrium points},
Int. J. Bifur. Chaos Appl. Sci. Eng. 30, 2050157, 8 pp (2020).

\bibitem{Fre} Freire, E., Ponce, E., Rodrigo, F. and Torres, F.: {\it Bifurcation sets of continuous piecewise linear systems with two zones}, Int. J. Bifur. Chaos Appl. Sci. Eng. 8, 2073-2097 (1998).

\bibitem{Fre1} Freire, E., Ponce, E., Rodrigo, F. and Torres, F.: {\it Bifurcation sets of symmetrical continuous piecewise linear systems with three zones}, Int. J. Bifur. Chaos Appl. Sci. Eng. 12, 1675-1702 (2002).

\bibitem{Fre3} Freire, E., Ponce, E. and Torres, F.: {\it A general mechanism to generate three limit cycles in planar Filippov systems with two zones}, Nonlin. Dyn. 78, 251-263 (2014).

\bibitem{Fre2} Freire, E., Ponce E. and Torres, F.: {\it Canonical discontinuous planar piecewise linear systems}, SIAM J. Appl. Dyn. Syst. 11, 181-211 (2012).

\bibitem{Fre4} Freire, E., Ponce, E. and Torres, F.: {\it The discontinuous matching of two planar linear foci can have three nested crossing limit cycles}, Publ. Mat. 2014, 221-253 (2014).

\bibitem{P5} Herssens, C., Maesschalck, P., Arte's, J. C., Dumortier, F., and Llibre, J.: {\it Piecewise polynomial planar phase portraits-P5}, \linebreak https://www.uhasselt.be/uh/dysy/software/p5.html, 2008.

\bibitem{Hua3} Huan, S. and Yang, X.: {\it Existence of limit cycles in general planar piecewise linear systems of saddle-saddle dynamics}, Nonlin. Anal. 92, 82-95 (2013).

\bibitem{Hua2} Huan, S. and Yang, X.: {\it Generalized Hopf bifurcation emerged from a corner in general planar piecewise smooth systems}, Nonlin. Anal. 75, 6260-6274 (2012).

\bibitem{Hua} Huan, S. and Yang, X.: {\it On the number of limit cycles in general planar piecewise linear systems of node-node types}, J. Math. Anal. Appl. 411, 340-353 (2014).

\bibitem {Li} Li, L.: {\it Three crossing limit cycles in planar piecewise linear systems with saddle focus type}, Electron. J. Qual. Theory Differ. Equ. 70, 1-14 (2014).

\bibitem {Li2} Li, S. and Llibre, J.: {\it On the limit cycles of planar discontinuous piecewise linear differential systems with a unique equilibrium}, Discrete Contin. Dyn. Syst. B 24, 5885-5901 (2019).

\bibitem {Li3} Li, S. and Llibre, J.: {\it Phase portraits of piecewise linear continuous differential systems with two zones separated by a straight line}, J. Differ. Equ. 266, 8094-8109 (2019).

\bibitem{Lli4} Llibre, J., Ponce, E. and Valls, C.: {\it Uniqueness and non-uniqueness of limit cycles for piecewise linear differential systems with three zones and no symmetry}, J. Nonlin. Sci. 25, 861-887 (2015).

\bibitem{Lli1} Llibre, J. and Teixeira, M.: {\it Piecewise linear differential systems without equilibria produce limit cycles?}, Nonlin. Dyn. 88, 157-167 (2017).

\bibitem{Lli3} Llibre, J. and Teixeira, M.: {\it Piecewise linear differential systems with only centers can create limit cycles?}, Nonlin. Dyn. 91, 249-255 (2018).

\bibitem{Lli2} Llibre, J. and Zhang, X.: {\it Limit cycles for discontinuous planar piecewise linear differential systems separated by one straight line and having a center}, J. Math. Anal. Appl. 467, 537-549 (2018).

\bibitem{Chu} Lum, R. and Chua, L.: {\it Global properties of continuous piecewise linear vector fields. Part I: Simplest case in $\R^2$}, Memorandum UCB/ERL M90/22, University of California at Berkeley (1990).

\bibitem{Med} Medrado, J. and Torregrosa, J.: {\it Uniqueness of limit cycles for sewing planar piecewise linear systems}, J. Math. Anal. Appl. 431, 529-544 (2015). 

\bibitem{Mer} Mereu, A., Oliveira, R. and Rodrigues, C.: {\it Limit cycles for a class of discontinuous piecewise generalized Kukles differential systems}, Nonlin. Dyn. 93, 2201-2212 (2018).

\bibitem{Pon} Ponce, E., Ros, J. and Vela, E.: {\it Limit cycle and boundary equilibrium bifurcations in continuous planar piecewise linear systems}, Int. J. Bifur. Chaos Appl. Sci. Eng. 25, 1530008, 18 pp (2015).

\bibitem{Pu} Pu, J., Chen, X., Chen, H. and Xia, Y.: {\it Global analysis of an asymmetric xontinuous piecewise linear differential aystem with three linear zones},  Int. J. Bifur. Chaos Appl. Sci. Eng. 31, 2150027, 41 pp (2021).

\bibitem{Wan2} Wang, J., Chen, X. and Huang, L.: {\it The number and stability of limit cycles for planar piecewise linear systems of node-saddle type}, J. Math. Anal. Appl. 469, 405-427 (2019).

\bibitem{Wan1} Wang, J., He, S. and Huang, L.: {\it Limit cycles induced by threshold nonlinearity in planar piecewise linear systems of node-focus or node-center type}, Int. J. Bifur. Chaos Appl. Sci. Eng. 30, 2050160, 13 pp (2020).

\bibitem{Wan3} Wang, J., Huang, C. and Huang, L.: {\it Discontinuity induced limit cycles in a general planar piecewise linear system of saddle-focus type}, Nonlin. Anal. Hybrid Syst. 33, 162-178 (2019).

\bibitem{Mathematica} Wolfram Research, Inc., \emph{Mathematica}, Version 12.2, Champaign, IL (2020). https://www.wolfram.com/mathematica

\bibitem{Xio2} Xiong, Y. and Han, M.: {\it Limit cycle bifurcations in discontinuous planar systems with multiple lines}, J. Appl. Anal. Comput. 10, 361-377 (2020).

\bibitem{Xio} Xiong, Y. and Wang, C.: {\it Limit cycle bifurcations of planar piecewise differential systems with three zones}, Nonlin. Anal. Real World Appl. 61, 103333, 18 pp (2021).
 
 
\end{thebibliography}

\end{document}